\documentclass[a4paper]{article}
\title{$P$-adic $L$-functions of Bianchi modular forms}
\date{}
\author{Chris Williams}

\usepackage[margin=1.29in]{geometry}
\usepackage{preamble}

\begin{document}

\maketitle

%
%
\begin{abstract}The theory of overconvergent modular symbols, developed by Rob Pollack and Glenn Stevens, gives a beautiful and effective construction of the $p$-adic $L$-function of a  modular form. In this paper, we give an analogue of their results for Bianchi modular forms, that is, modular forms over imaginary quadratic fields. In particular, we prove control theorems that say that the canonical specialisation map from overconvergent to classical Bianchi modular symbols is an isomorphism on small slope eigenspaces of suitable Hecke operators. We also give an explicit link between the classical modular symbol attached to a Bianchi modular form and critical values of its $L$-function, which then allows us to construct $p$-adic $L$-functions of Bianchi modular forms.
\end{abstract}

%
%
\section*{Introduction}
In recent years, $p$-adic $L$-functions have been a subject of considerable study. As an example, the `main conjecture' of Iwasawa theory says that $p$-adic $L$-functions control the size of cohomology groups of Galois representations. They have also been used in the construction of Stark--Heegner points, which are conjecturally-global points on non-CM elliptic curves over number fields. Via the theory of overconvergent modular symbols, Rob Pollack and Glenn Stevens (see \cite{PS11} and \cite{PS12}, or for an exposition, \cite{Pol11}) gave a beautiful and computationally effective method of constructing the $p$-adic $L$-function of a suitable rational modular form. They defined the space of \emph{overconvergent modular symbols} to be a space of modular symbols taking values in some $p$-adic distribution space $\D_k(L)$, exhibiting a surjective Hecke-equivariant specialisation map from this space to the space of classical modular symbols of weight $k$; the crux of their work is Stevens' control theorem, which says that this map is an isomorphism on the small slope eigenspaces of the $U_p$-operator. This draws comparison with the work of Coleman on overconvergent modular forms, and in particular his result that a small slope overconvergent modular form is in fact classical.\\
\\
In this paper, we give an analogue of their theory for the case of Bianchi modular forms, that is, the case of modular forms over imaginary quadratic fields. The theory of modular symbols over the rationals generalises to arbitrary number fields; however, the definitions can be made far more explicit over $\Q$ or imaginary quadratic fields, motivating the study of the imaginary quadratic case as a stand alone topic. We give a brief description of the completely general case as motivation for the sequel. To this end, let $F$ be a number field with adele ring $\A_F$ (with finite adeles $\A_F^f$) and degree $d_F = r + 2s$, where $r$ (resp. $s$) is the number of real (resp. pairs of complex) embeddings of $F$. Let $\Omega$ be an open compact subgroup of $\GLt(\A_F^f)$, which takes the (adelic) role of a congruence subgroup. The analogue of the modular curve is then the locally symmetric space
\[Y_\Omega \defeq \GLt(F)\backslash\GLt(\A_F)/\Omega K_\infty,\]
where $K_\infty$ is the standard maximal compact subgroup of $\GLt(F\otimes_{\Q}\R)$. Let $q = r+s$. For a suitable local system $\mathcal{V}$ on $Y_\Omega$, we have a canonical link between spaces of automorphic forms of fixed weight and level $\Omega$ and elements of the cohomology group $\h_{\mathrm{c}}^q(Y_\Omega,\mathcal{V})$. We call this latter space the space of \emph{modular symbols of level $\Omega$ and weight $\mathcal{V}$}. When $\det(\Omega) = (\roi_K\otimes \widehat{\Z})^\times,$ and via a result called strong approximation, $Y_\Omega$ breaks down into components indexed by the class group of $F$; explicitly, we choose representatives for the class group, which then allows us to define a collection $(\Gamma_\alpha)_{\alpha \in \mathrm{Cl}(F)}$ of discrete subgroups of $\SLt(F)$, and a collection of spaces $(Y_\alpha)_{\alpha\in\mathrm{Cl}(F)}$, where $Y_\alpha \defeq \Gamma_\alpha\backslash[\uhp^r\times\uhs^s]$. We see that $Y_\Omega = \bigsqcup_{\alpha \in \cl(F)}Y_{\alpha}$, and this leads to a decomposition
\[\h_{\mathrm{c}}^q(Y_\Omega,\mathcal{V}) \cong \bigoplus_{\alpha \in \mathrm{Cl}(F)}\h_{\mathrm{c}}^q(Y_\alpha,\mathcal{V}).\]
In the Bianchi case, where $q = 1$, we can study the geometry of $\uhs$ to find that each of these components is isomorphic to a space with a much simpler description, namely the natural analogue of that used in \cite{PS11}, where they define modular symbols to be functions on divisors on cusps into some polynomial space. This isomorphism -- described in \cite{AS86} -- only exists for $q=1$, showing that this more explicit description of modular symbols only applies in the rational and imaginary quadratic cases. In this paper, we give a complete account of this explicit description in the latter setting.\\
\\
The literature regarding Bianchi modular forms is widespread; in particular, an account of the general theory over arbitrary number fields is given in Andr\'{e} Weil's book \cite{Wei71}, whilst accounts in the imaginary quadratic case for weight 2 are given by John Cremona and two of his students in \cite{Cre81}, \cite{CrWh94} and \cite{Byg98}. However, there seem to be a large number of different conventions pertaining to the theory surrounding these objects, and as such, we start with a brief introduction to the theory -- largely following the account of Eknath Ghate in \cite{Gha99} -- fixing as we do so the conventions that we'll use in the sequel. Broadly speaking, a \emph{Bianchi modular form} is an adelic automorphic form for $\GLt$ over an imaginary quadratic field $K$; such a form corresponds to a collection of automorphic functions $F^1, ..., F^h:\GLt(\C) \rightarrow V_{2k+2}(\C)$, where $h$ is the class number of $K$ and $V_k(\C)$ is the space of homogeneous polynomials in two variables of degree $k$ over $\C$, and then further to functions $\f^1,...,\f^h$ on the \emph{upper half-space} $\uhs \defeq \C\times\R_{>0}$, the analogue of the upper half-plane in this setting. Such a Bianchi modular form has a Fourier expansion, and accordingly an $L$-function that converges absolutely on some right half-plane. As in the case of classical modular forms, one can give a formula for the L-function in terms of integrals of the automorphic form over certain paths in $\uhs$. In Theorem \ref{integralformula}, we prove a precise version of:
\begin{mthmnum}\label{integralformula}Let $\Phi$ be a cuspidal Bianchi modular form of weight $(k,k)$ and level $\Omega_1(\n)$. Let $\psi$ be a Grossencharacter with a prescribed infinity type. Then there is an integral formula for $L(\Phi,\psi,s)$ for $s\in \C$.
\end{mthmnum}
The heart of the paper is in the study of Bianchi modular symbols. To a cuspidal Bianchi modular form $\Phi$ of weight $(k,k)$ and level $\Omega_1(\n)$, where $(p)|\n$, we associate a collection $\f^1, ..., \f^h$ of functions on $\uhs$, as above, each satisfying an automorphy condition for some discrete subgroup $\Gamma_i$ of $\SLt(K)$. To each of these $\f^i$, we associate a classical $V_k(\C)\otimes_{\C}V_k(\C)$-valued modular symbol $\phi_{\f^i}$ for $\Gamma_i$, and exhibit a link between values of this symbol and critical values of the part of the $L$-function corresponding to $\f^i$. In particular, we have the following (see Theorem \ref{lfunctionmodsymb} for a precise statement):
\begin{mthmnum}Let $\Phi$ be a cuspidal Bianchi modular form of parallel weight $(k,k)$ and level $\Omega_1(\n)$ with renormalised $L$-function $\Lambda(\Phi,\psi)$. For a Grossencharacter $\psi$ of $K$ of conductor $\ff$ and infinity type $0 \leq (q,r) \leq (k,k)$, we can express $\Lambda(\Phi,\psi)$ explicitly in terms of the coefficients of $\phi_{\f^i}(\{\alpha\}-\{\infty\}),$ for $1\leq i \leq h$ and $\alpha$ ranging over a set of of representatives of $\ff^{-1}/\roi_K$ with $\alpha\ff$ coprime to $\ff$. Here the coefficients are taken in a suitable basis of $V_k(\C)\otimes_{\C}V_k(\C)$.
\end{mthmnum}
We then define the space of overconvergent Bianchi modular symbols to be the space of modular symbols taking values in some $p$-adic distribution space; precisely, we fix a finite extension $L/\Qp$, and denoting by $\A_k(L)$ the space of rigid analytic functions on the unit disc defined over $L$, our distribution space is $\D_{k,k}(L)\defeq \text{Hom}(\A_k(L)\ctp_{L} \A_k(L),L)$. After renormalising the values of $\phi_{\f^i}$, we may consider it to have values in $V_k(L)\otimes_{L}V_k(L)$, and we then have a specialisation map from overconvergent to classical modular symbols by dualising the inclusion $V_k(L)\otimes_{L}V_k(L) \hookrightarrow \A_k(L)\ctp_{L}\A_k(L)$, much like in the rational case. We prove the following analogue of Stevens' control theorem; in the case of $p$ inert, it is proved in Corollary \ref{contthm1}, in the case of $p$ ramified it is proved in Corollary \ref{contthm1} combined with Lemma \ref{untou}, and in the case $p$ split it is proved in Theorem \ref{refinedcontrolthm}.
\begin{mthmnum}[The Control Theorem for Bianchi Modular Symbols]Let $p$ be a rational prime with $p\roi_K = \prod_{\pri|p}\pri^{e_{\pri}}$. For each prime $\pri|p$, let $\lambda_{\pri} \in L$. Then, when $v_p(\lambda_{\pri}) < (k+1)/e_{\pri}$ for all $\pri|p$, the restriction of the map 
\[\rho:\bigoplus_{i=1}^h\symb_{\Gamma_i}(\D_{k,k}(L))^{\{U_{\pri}=\lambda_{\pri} : \pri|p\}} \longrightarrow \bigoplus_{i=1}^h\symb_{\Gamma_i}(V_{k,k}^*(L))^{\{U_{\pri}=\lambda_{\pri} : \pri|p\}}\]
to the simultaneous $\lambda_{\pri}$-eigenspaces of the $U_{\pri}$-operators is an isomorphism.
\end{mthmnum}
 The proof draws from work of Matthew Greenberg in \cite{Gre07}, in that we define a series of \emph{finite approximation modules}, and lift compatibly through this system to obtain a overconvergent symbol from a classical one.\\
\\
It is worth remarking that whilst in the rational case, the control theorem gives an analogue of Coleman's small slope classicality theorem, no such theory of `overconvergent Bianchi modular forms' yet exists.\\
\\
In the remainder, the values of an overconvergent eigenlift $\Psi_{\f^i}$ are studied; namely, we prove that such a symbol takes values in some space of locally analytic distributions, and that it is \emph{admissible} (or \emph{tempered}).\\
\\
Suppose we are in the set-up of the control theorem, and let $\Phi$ be a small slope cuspidal Bianchi eigenform of weight $(k,k)$ and level $\Omega_1(\n)$. Then to $\Phi$ we can associate a small slope eigensymbol $(\phi_1,...,\phi_h)$ in a direct sum of symbol spaces, which we can lift uniquely to an overconvergent symbol $(\Psi_1,...,\Psi_h)$ using the control theorem. Then there is a way of patching together the distributions $\Psi_i(\{0\}-\{\infty\})$ to a locally analytic distribution $\mu_p$ on the ray class group $\cl(K,p^\infty)$. We then have the following (see Theorem \ref{padiclfn}):

\begin{mthmnum}
Let $\Phi$ be a cuspidal Bianchi eigenform of weight $(k,k)$ and level $\Omega_1(\n)$, where $(p)|\n$, with $U_\pri$-eigenvalues $a_\pri$, where $v(a_{\pri})<(k+1)/e_{\pri}$ for all $\pri|p$. Then there exists a locally analytic distribution $\mu_p$ on $\cl(K, p^\infty)$ such that for any Grossencharacter of conductor $\ff|(p^\infty)$ and infinity type $0 \leq (q,r) \leq (k,k)$, we have
\[\mu_p(\psi_{p-\mathrm{fin}}) = A(\psi)\Lambda(\Phi,\psi),\]
where $\psi_{p-\mathrm{fin}}$ is a character on $\cl(K,p^\infty)$ corresponding to $\psi$, $A$ is a scalar depending only on $\psi$ and given explicitly in Theorem \ref{padiclfn}, and $\Lambda(\Phi,\psi)$ is the normalised $L$-function of $\Phi$. The distribution $\mu_p$ is $(h_{\pri})_{\pri|p}$-admissible, where $h_{\pri} = v_p(a_{\pri})$, and hence is unique.
\end{mthmnum}
We call $\mu_p$ the \emph{$p$-adic $L$-function of $\Phi$}. Such a result is the natural analogue of the results of Pollack and Stevens in the rational case.\\
\\
\textbf{Comparison to relevant literature:} There are a number of people who have worked on similar things in the recent past. Perhaps of most relevance is Mak Trifkovic, who in \cite{Tri06} performed computations with overconvergent Bianchi modular symbols. He proved a lifting theorem in the case of weight 2 ordinary eigenforms over an imaginary quadratic field of class number 1, using similar explicit methods to \cite{Gre07}. The lifting results in this paper are a significant generalisation of his theorem, though the author has not made any efforts to repeat the computational aspects of Trifkovic's work in this more general setting. Looking in a different direction, in his PhD thesis (\cite{Bar13}) Daniel Barrera Salazar has generalised the results of Pollack and Stevens to the case of Hilbert modular forms (that is, to modular symbols over totally real fields), though the methods he uses are vastly different to those in the Bianchi case, as such symbols live in higher compactly supported cohomology groups and don't lend themselves to such explicit study.\\
\\
\textbf{Acknowledgments:} I would like to thank my PhD supervisor David Loeffler for suggesting this topic to me, as well as for the many conversations we've had on the subject. Thanks also to John Cremona and Haluk Sengun for helpful conversations relating to the classical theory of Bianchi modular forms.

\section{Bianchi Modular Forms}\label{BMF}
A \emph{Bianchi modular form} is an automorphic form for $\GLt$ over an imaginary quadratic field. Here, we give only a very brief description of the theory; more detailed descriptions are given by Jeremy Bygott (\cite{Byg98}, focusing on weight 2) and Eknath Ghate (\cite{Gha99}, for higher weights). The basic definitions are given in Section 1.1, whilst Section 1.2 looks at the $L$-function of a Bianchi modular form, giving an integral formula for the twisted $L$-function.
\begin{mnot}
Throughout this paper, we'll take $K$ to be an imaginary quadratic field with ring of integers $\roi_K$, different $\mathcal{D}$ and discriminant $-D$, $p$ a rational prime, $\n$ an ideal of $\roi_K$ divisible by $(p)$. At each prime $\pri$ of $K$, denote by $K_{\pri}$ the completion of $K$ with respect to $\pri$ and $\roi_{\pri}$ the ring of integers of $K_{\pri}$. Denote the adele ring of $K$ by $\A_K = K_\infty \times \A_K^f$, where $K_\infty$ are the infinite adeles and $\A_K^f$ are the finite adeles. Furthermore, define $\widehat{\roi_K} \defeq \roi_K \otimes_{\Z}\widehat{\Z}$ to be the finite integral adeles. Denote the class group of $K$ by $\cl(K)$ and the class number of $K$ by $h$, and -- once and for all -- fix a set of representatives $I_1, ..., I_h$ for $\cl(K)$, with $I_1 = \roi_K$ and each $I_j$ for $2\leq j\leq h$ integral and prime, with each $I_i$ coprime to $\n$ and $\mathcal{D}$. (We will assume that for any Grossencharacter of conductor $\ff$ considered in the sequel, these representatives are also coprime to $\ff$).\\
\\
Let $k \geq -1$ be an integer, and for any ring $R$, let $V_{2k+2}(R)$ denote the ring of homogeneous polynomials over $R$ in two variables of degree $2k+2$. Note that $V_{2k+2}(\C)$ is an irreducible complex right representation of $\SUt(\C)$, and denote the corresponding antihomomorphism by $\rho':\SUt(\C) \rightarrow \mathrm{GL}(V_{2k+2}(\C))$. Finally, define an antihomomorphism $\rho : \SUt(\C) \times \C^\times \rightarrow \mathrm{GL}(V_{2k+2}(\C))$ by $\rho(u,z) = \rho'(u)|z|^{-k}$.\\
\\
For a general Grossencharacter $\psi$ of $K$, for each prime $\pri$ of $K$ we denote by $\psi_{\pri}$ the restriction of $\psi$ to $K_{\pri}^\times$; then we also write $\psi_\infty$ for the restriction of $\psi$ to the infinite ideles, and $\psi_f$ for the restriction to the finite ideles.
\end{mnot}

\subsection{Definitions}
In the spirit of Weil, we define Bianchi modular forms adelically. First, we fix a level.
\begin{mdef}\label{omega}
Define
\[\Omega \defeq \Omega_1(\n) \defeq \left\{\matrd{a}{b}{c}{d} \in \GLt(\widehat{\roi_K}): c \in \n\widehat{\roi_K}, a,d \in 1 + \n\widehat{\roi_K}\right\}.\]
\end{mdef}
\begin{mdef}We say a function $\Phi:\GLt(\A_K) \rightarrow V_{2k+2}(\C)$ is a \emph{cusp form of weight $(k,k)$ and level $\Omega_1(\n)$} if it satisfies:
\begin{itemize}
\item[(i)]$\Phi(zgu) = \Phi(g)\rho(u,z)$ for $u \in \SUt(\C)$ and $z \in Z(\GLt(\C)) \cong \C^\times,$
\item[(ii)]$\Phi$ is right-invariant under the group $\Omega_1(\n)$,
\item[(iii)] $\Phi$ is left-invariant under $\GLt(K)$,
\item[(iv)] $\Phi$ is an eigenfunction of the operator $\partial$,
\[\partial f= (k^2/2 + k)f,\]
where $\partial/4$ denotes a component of the Casimir operator in the Lie algebra $\mathfrak{s}\mathfrak{l}_2(\C)\otimes_{\R}\C$ (see \cite{Hid93}, section 1.3), and where we consider $\Phi(g_\infty g_f)$ as a function of $g_\infty \in \GLt(\C)$, and
\item[(v)] $\Phi$ satisfies the cuspidal condition that for all $g \in \GLt(\A_K)$, we have
 \[\int_{K\backslash\A_K}\Phi(ug)du = 0,\]
where we consider $\A_k$ to be embedded inside $\GLt(\A_K)$ by the map sending $u$ to $\smallmatrd{1}{u}{0}{1}$, and $du$ is the Lebesgue measure on $\A_K$.
\end{itemize}
\end{mdef}
\begin{mrems}
\begin{itemize}
\item[(i)] The cuspidal condition is a natural one; the value of the integral for a fixed $g$ corresponds to a constant Fourier coefficient.
\item[(ii)] The general definitions given above are already slightly tailored to work with cusp forms, that appear only at parallel weights (see \cite{Har87}); more generally, it is possible to define general automorphic forms of weight $(k_1,k_2)$, for distinct integers $k_1$ and $k_2$, over $K$. 
\end{itemize}
\end{mrems}

A cusp form $\Phi$ of weight $(k,k)$ and level $\Omega_1(\n)$ descends to give a collection of $h$ functions $F^i:\GLt(\C) \rightarrow V_{2k+2}(\C)$, for $i = 1,...,h$, where $h$ is the class number. Such a collection is non-canonical, depending on choices of representatives for the class group of $K$. To obtain this decomposition, we note that $\GLt(\A_K)$ decomposes as the disjoint union of $h$ sets. To describe this, recall that we took $I_1, ... ,I_h$ to be a complete set of (prime or trivial) representatives for the class group; Set $t_1 = 1$ and for each $i\geq 2$, define
\[t_i = (1,...,1,\pi_i,1,...) \in \A_K^\times,\]
where $\pi_i$ is a uniformiser in $K_{I_i}$. Then define
\[g_i = \matrd{1}{0}{0}{t_i} \in \GLt(\A_K).\]
Then we have:
\begin{mthm}[Strong Approximation]\label{strongapprox} There is a decomposition
\[\GLt(\A_K) = \coprod\limits_{i=1}^h \GLt(K)\cdot g_i\cdot \left[\GLt(\C)\times \Omega_1(\n)\right].\]
\end{mthm}
\begin{proof}
See \cite{Byg98}, Section 5.2.
\end{proof}
It is now clear that $\Phi$ descends in the way claimed above via $F^i(g) \defeq \Phi(g_i g)$. Furthermore, we can descend further using condition (i) to obtain $h$ functions $\f^i : \uhs \rightarrow V_{2k+2}(\C)$, where
\[\uhs \defeq \GLt(\C)/\left[\SUt(\C)\cdot Z(\GLt(\C))\right]\]
is hyperbolic 3-space. Such functions are \emph{automorphic forms on $\uhs$ of weight $(k,k)$ and level }
\begin{align}\label{gammai}\Gamma_i \defeq \SLt(K) \cap g_i\Omega_1(\n)g_i^{-1}.
\end{align}
This process is described in \cite{Byg98}. Note here that $\Gamma_1$ is nothing other than the usual $\Gamma_1(\n) \leq \SLt(\roi_K)$.

\subsection{$L$-functions}
\subsubsection{Fourier Expansions}\label{fourierexpansions}
Let $\Phi:\GLt(\A_K)\rightarrow V_{2k+2}(\C)$ be a cusp form of weight $(k,k)$ and level $\Omega_1(\n)$. Then, from \cite{Hid94}, $\Phi$ has a Fourier expansion of the form 
\begin{align}
\label{fourierexpansionad}
\Phi\left[\matrd{t}{z}{0}{1}\right] = |t|_K \sum\limits_{\alpha \in K^\times} c(\alpha t\delta, \Phi)W(\alpha t_\infty)e_K(\alpha z),
\end{align}
where:
\begin{itemize}
\item[(i)]$|\cdot|_K$ is the idelic norm character,
\item[(ii)]$\delta = \sqrt{-D}$ (where $-D$ is the discriminant of $K$) is a generator of the different $\mathcal{D}$ of $K$, i.e. $\delta\roi_K = \mathcal{D}$,
\item[(iii)]The Fourier coefficient $c(\cdot, \Phi)$ may be considered as a function on the fractional ideals of $K$, with $c(I,\Phi) = 0$ for $I$ non-integral,
\item[(iv)]$e_K$ is an additive character of $K\backslash\A_K$ defined by
\[e_K = \left(\prod\limits_{\pri\text{ prime}}(e_p\circ\mathrm{Tr}_{K_{\pri}/\Q_p})\right)\cdot(e_\infty\circ\mathrm{Tr}_{\C/\R}),\]
for
\[e_p\left(\sum\limits_{j}d_jp^j\right) = e^{-2\pi i\sum\limits_{j<0}d_jp^j} \hsp\text{ and }\hsp e_\infty(r) = e^{2\pi ir},\]
and
\item[(v)]$W:\C^\times \rightarrow V_{2k+2}(\C)$ is the \emph{Whittaker function}
\[W(s) \defeq \sum\limits_{n=0}^{2k+2}\binomc{2k+2}{n}\left(\frac{s}{i|s|}\right)^{k+1-n}K_{n-(k+1)}(4\pi |s|)X^{2k+2-n}Y^n,\]
where $K_n(x)$ is a (modified Bessel function that is) a solution to the differential equation
\[\frac{d^2K_n}{dx^2} + \frac{1}{x}\frac{dK_n}{dx} - \left(1+ \frac{n^2}{x^2}\right)K_n = 0,\]
with asymptotic behaviour
\[K_n(x) \sim \sqrt{\frac{\pi}{2x}}e^{-x}\]
as $x \rightarrow \infty.$ Note that, in particular, $K_{-n} = K_n$.
\end{itemize}

If our adelic cusp form $\Phi$ corresponds to a collection of $h$ cusp forms $\f^1, ..., \f^h$ on $\uhs$, where $h$ is the class number of $K$, then the Fourier expansion stated above descends to give a Fourier expansion for each $\f^j$, which can be worked out to be (see \cite{Gha99})
\begin{align}\label{fourierexpansion}
\f^j\left(z,t;\binomc{X}{Y}\right) &= |t_j|_Kt\sum\limits_{n=0}^{2k+2}\binomc{2k+2}{n}\f^j_n(z,t)X^{2k+2-n}Y^n,\\ 
\f^j_n(z,t) &\defeq \sum\limits_{\alpha\in K^\times}\bigg[c(\alpha t_j\delta)\left(\frac{\alpha}{i|\alpha|}\right)^{k+1-n}K_{n-k-1}(4\pi|\alpha|t)e^{2\pi i (\alpha z + \overline{\alpha z})}\bigg].\notag\end{align}

\subsubsection{Defining the $L$-Function}\label{deflfn}
\begin{mdef}
Let $\Phi$ be a cusp form for $K$ of any weight and level $\Omega_1(\n)$, with Fourier expansion given by equation (\ref{fourierexpansionad}). The \emph{$L$-function} of $\Phi$ is defined by
\[L(\Phi,s) \defeq \sum\limits_{0\neq \m \subset \roi_K}c(\m,\Phi)N(\m)^{-s},\]
where the sum is over all non-zero ideals of $\roi_K$.
\end{mdef}
Let $\f^1, ..., \f^h$ be the $h$ automorphic forms corresponding to $\Phi$ with respect to $I_1,...,I_h$. Then to each one, we assign a `part' of this $L$-function, in the following sense:
\begin{mdef}
Let $w = |\roi_K^\times|$. Define 
\[L^i(\Phi,s) = L(\f^i,s) \defeq w^{-1} \sum\limits_{\alpha\in K^\times}c(I_i(\alpha), \Phi)N(I_i(\alpha))^{-s}.\]
\end{mdef}
Note here that
\[L(\Phi,s) = L^1(\Phi,s) + \cdots + L^h(\Phi,s),\]
where here we scale by $w^{-1}$ as when we sum over elements of $K^\times$, we include each ideal $w$ times (once for each unit).\\
\\
The following lemma (see \cite{Wei71}, Chapter II) shows that each component of the $L$-function is a holomorphic function on a right half-plane:
\begin{mlem}\label{absconv}
For each $i$, $L^i(\Phi,s)$ converges absolutely in the right half-plane Re$(s) > C$, for some constant $C$.
\end{mlem}
Let $\psi$ be a Grossencharacter of conductor $\ff$. We want to define the twist of the $L$-function by $\psi$. For each ideal $\m$ coprime to $\ff$, define $\psi(\m) \defeq \prod_{\pri^n||\m}\psi_{\pri}(\pi_{\pri})^n$; then define
\[L(\Phi,\psi,s) \defeq \sum\limits_{\substack{0\neq \m \subset \roi_K\\ (\m,\ff) = 1}}c(\m,\Phi)\psi(\m)N(\m)^{-s},\]
where, for ideals $I$ and $J$, the notation $(I,J) = 1$ means that $I$ and $J$ are coprime. It is more convenient in terms of our future work to define the twist in terms of the $L^i(\Phi,s).$ Recall that $\mathcal{D} = (\delta)$ is the different of $K$ and the definition of $L^i(\Phi,s)$ above, and set
\begin{align*}L^i(\Phi,\psi,s)&= w^{-1}\sum\limits_{\alpha\in K^\times}c(\alpha\delta I_i,\Phi)\psi(\alpha\delta I_i)N(\alpha\delta I_i)^{-s}\\ 
&= w^{-1}\psi(t_i)|t_i|_f^s\sum\limits_{\substack{\alpha\in \mathcal{D}^{-1}\\(\alpha\delta I_i,\ff) = 1}}c(\alpha\delta I_i,\Phi)\psi((\alpha\delta))|\alpha\delta|^{-2s}
\end{align*}
where here $|\alpha\delta|$ is the usual norm of $\alpha\delta \in K^\times$ and $|t_i|_f = \prod_{v\nmid\infty}|t_i|_v$ is the finite idelic norm.
\begin{mrem}
Note that in \cite{Wei71}, Chapter II, it is proved that the twisted $L$-function still converges absolutely in some suitable right half-plane.
\end{mrem}
There is one more useful alternative method of writing the twist. Define
\[\psi_{\ff} = \prod_{\pri|\ff}\psi_{\pri}.\]
Then note that as $\psi$ is trivial on the diagonal embedding of $K^\times$ in $\A_K^\times$, we have, for $(\alpha)$ coprime to $\ff$,
\[\psi_\infty(\alpha)\psi_{\ff}(\alpha)\psi((\alpha)) = 1.\]
Accordingly, we can also write the twist as
\[L^i(\Phi,\psi,s) = w^{-1}\psi(t_i)|t_i|_f^s\sum\limits_{\substack{\alpha \in K^\times\\ (\alpha\delta I_i,\ff) = 1}}c(\alpha\delta I_i,\Phi)\psi_\infty(\alpha\delta)^{-1}\psi_{\ff}(\alpha\delta)^{-1}|\alpha\delta|^{-2s}.\]

\subsubsection{An Integral Formula}
\label{lfunction}
In the rational case, the twist of an $L$-function attached to a modular form can be written in terms of an integral formula, which is crucial to its study via modular symbols. We now obtain such an integral formula for $L^i(\Phi,\psi,s)$. To do so, we'll need some results about Gauss sums. Hecke wrote down a suitable analogue of Gauss sums over general number fields in \cite{Hec20} and \cite{Hec23}, and Deligne wrote down local analogues (via the theory of local $\varepsilon$-factors) in \cite{Del72}, Section 3; his work is translated into English in \cite{Tat79}, Section 3. In \cite{Nar04}, Proposition 6.14(ii), it is shown that Hecke's Gauss sums satisfy a product law, and moreover that the local factors agree with Deligne's local $\varepsilon$-factors. The version we'll use is one obtained from Hecke's definitions by Nemchenok in \cite{Nem93} Section 6. For a Grossencharacter $\psi$ of conductor $\ff$, the Gauss sum is defined to be
\[\tau(\psi) \defeq \sum\limits_{\substack{[a]\in\ff^{-1}/\roi_K\\ ((\alpha)\ff, \ff) = 1}}\psi(a\ff)\psi_\infty\left(\frac{a}{\delta}\right)e^{2\pi i\mathrm{Tr}_{K/\Q}(a/\delta)},\]
where as before $\delta$ is a generator of the different $\DD$ of $K$. The key proposition we require is the following, which follows easily from \cite{Nar04}, Proposition 6.14(i), and allows us to deal with the terms of the $L$-series which correspond to ideals not coprime to the conductor. 
\begin{mprop}\label{gausssum}
\begin{itemize}
\item[(i)] For all $b\in\roi_K$, we have
\[ \sum\limits_{\substack{[a]\in\ff^{-1}/\roi_K\\ ((a)\ff, \ff) = 1}}\psi(a\ff)\psi_\infty\left(\frac{a}{\delta}\right)e^{2\pi i\mathrm{Tr}_{K/\Q}(ab/\delta)} =  \tau(\psi)\psi_{\ff}(b).\]
\item[(ii)] By replacing $\psi$ with $\psi^{-1}$, we have
\[\frac{1}{\tau(\psi^{-1})}\sum\limits_{\substack{[a]\in\ff^{-1}/\roi_K\\ ((a)\ff, \ff) = 1}}\psi(a\ff)^{-1}\psi_\infty\left(\frac{a}{\delta}\right)^{-1}e^{2\pi i\mathrm{Tr}_{K/\Q}(ab/\delta)} = \left\{\begin{array}{ll}\psi_{\ff}(b)^{-1} & : ((b),\ff) = 1,\\
0 & : \text{otherwise}.
\end{array}\right.\]
\end{itemize}
\end{mprop}
From hereonin, $\psi$ denotes a Grossencharacter of infinity type $(-u,-v)$ and conductor $\ff$. Ultimately, the result we will prove is:
\begin{mthm}\label{integralformula}Let $\Phi$ be a cuspidal Bianchi modular form of weight $(k,k)$ and level $\Omega_1(\n)$, corresponding to $h$ cusp forms $\f^1,...,\f^h$ on $\uhs$. For $n \in \{0,...,2k+2\}$, let $\f^i_n$ be as defined in equation (\ref{fourierexpansion}) above. Let $\psi$ be a Grossencharacter with infinity type $(-u,-v) = (-\tfrac{k+1-n}{2},\tfrac{k+1-n}{2})$. Then, for $s\in \C$, we have
\[L^i(\Phi,\psi,s)=A(i,n,\psi,s)\sum\limits_{[a]\in\ff^{-1}/\roi_K}\psi(a\ff)^{-1}a^u\overline{a^v}\int_0^\infty t^{2s-2}\f^i_n(a,t)dt,\]
where
\[A(i,n,\psi,s) = \psi(t_i)|t_i|_f^{s-1}\cdot\frac{4\cdot(2\pi)^{2s} i^{k+1-n}\smallbinomc{2k+2}{n}^{-1}}{|\delta|^{2s}\itGamma\left(s+\tfrac{n-k-1}{2}\right)\itGamma\left(s - \tfrac{n-k-1}{2}\right)w\tau(\psi^{-1})}\]
is an explicit function of $s$.
\end{mthm}
The proof is contained in the rest of this section.\\
\\
Using Proposition \ref{gausssum} to replace $\psi_{\ff}(\alpha\delta)^{-1}$ in the expression for $L^i(\Phi,\psi,s)$ we obtained above, and rearranging, we get
\[L^i(\Phi,\psi,s) = \frac{\psi(t_i)|t_i|_f^s}{w\tau(\psi^{-1})|\delta|^{2s}}\sum\limits_{\substack{[a]\in\ff^{-1}/\roi_K\\ ((a)\ff, \ff) = 1}}\psi(a\ff)^{-1}a^u\overline{a^v}
\sum\limits_{\alpha\in K^\times}
c(\alpha\delta I_i,\Phi)\alpha^u\overline{\alpha^v} e^{2\pi i\mathrm{Tr}_{K/\Q}(a\alpha)}|\alpha|^{-2s}.
\]

To get our integral formula, we will use the standard integral (see \cite{Hid94}, section 7)
\[\int_0^\infty t^{j-1}K_{n-k-1}(\lambda t)dt = \lambda^{-j}2^{j-2}\mathit{\Gamma}\left(\frac{j+n-k-1}{2}\right)\mathit{\Gamma}\left(\frac{j-n+k+1}{2}\right).\]
Setting $\lambda = 4\pi|\alpha|$, and $j = 2s$, we get an expression for $|\alpha|^{-2s}.$ Now, to force our expressions into the form of the Fourier expansion for $\f^i$, we fix $(u,v) = (\tfrac{k+1-n}{2},-\tfrac{k+1-n}{2})$; then substituting the expression for $|\alpha|^{-2s}$ and rearranging the terms, we get that this is equal to

\[=A(i,n,\psi,s)\sum\limits_{\substack{[a]\in\ff^{-1}/\roi_K\\ ((a)\ff, \ff) = 1}}\psi(a\ff)^{-1}a^u\overline{a^v}\int_0^\infty t^{2s-2}\f^i_n(a,t)dt,\]
where
\[A(i,n,\psi,s) = \frac{\psi(t_i)|t_i|_f^{s-1}}{|\delta|^{2s}}\cdot\frac{4\cdot(2\pi)^{2s} i^{k+1-n}\smallbinomc{2k+2}{n}^{-1}}{\itGamma\left(s+\tfrac{n-k-1}{2}\right)\itGamma\left(s - \tfrac{n-k-1}{2}\right)w\tau(\psi^{-1})}\]
is an explicit function of $s$.

\section{Bianchi Modular Symbols}
In this section, we describe Bianchi modular symbols. Recalling the general motivation given in the introduction, we elaborate in detail on the Bianchi case, focusing entirely on the description in terms of functions rather than compactly supported cohomology groups. We'll also show how to obtain a Bianchi modular symbol from a Bianchi modular form; the process is, unfortunately, significantly less intuitive than the rational case. We conclude the section by exhibiting a link between the modular symbol of a Bianchi modular form and critical values of its $L$-function.
\subsection{Definitions}
\label{modsymbolsection}
\begin{mdef}
We define, for a non-negative integer $k$ and a ring $R$, the space $V_k(R)$ as above to be the space of homogeneous polynomials of degree $k$ in two variables over $R$. Furthermore, we define $V_{k,k}(R) \defeq V_k(R) \otimes_{R}V_k(R).$
\end{mdef}
Note that we can identify $V_{k,k}(\C)$ with the space of polynomials that are homogeneous of degree $k$ in two variables $X,Y$ and homogeneous of degree $k$ in two further variables $\Xbar, \Ybar$. Furthermore, $V_k(\C)$ is an irreducible $\SUt(\C)$-module, with $\SUt(\C)$ acting on the right by
\[P|\gamma\binomc{X}{Y} = P\left(\gamma\binomc{X}{Y}\right).\]
The following defines a different action on this space, with a view to obtaining a `nice' action on the dual space of $V_{k,k}(\C)$.
\begin{mdef}\label{actionofsltc}We have a left-action of $\SLt(\C)$ on $V_k(\C)$ defined by
\[\gamma\cdot P\binomc{X}{Y} = P\binomc{dX +bY}{cX + aY}, \hsp \gamma = \matr.\]
We then obtain a left-action of $\SLt(\C)$ on $V_{k,k}(\C)$ by
\[\gamma\cdot P\left[\binomc{X}{Y},\binomc{\Xbar}{\Ybar}\right] = P\left[\binomc{dX + bY}{cX + aY},\binomc{\bar{d}\Xbar + \bar{b}\Ybar}{\bar{c}\Xbar + \bar{a}\Ybar}\right].\]
\end{mdef}
\begin{mrem}
This action, whilst appearing unconventional, is chosen so that it is compatible with an action on a space of locally analytic functions on some $p$-adic space; see Section \ref{actionofsig} for further details. This compatibility simplifies matters considerably when considering specialisation maps from overconvergent to classical modular symbols.
\end{mrem}
The left action of $\SLt(\C)$ described above translates into a \emph{right}-action on the dual space $V_{k,k}^*(\C) \defeq \Hom(V_{k,k}(\C),\C).$ For $\mu \in V_{k,k}^*(\C)$, we set
\[\mu|\gamma(P) = \mu(\gamma\cdot P).\]
We also have a left-action of $\SLt(\C)$ on the space $\Delta_0 \defeq \divzeroK$ by fractional linear transformations. 
\begin{mdef}
Let $\Gamma \leq \SLt(K)$ be a discrete subgroup, and let $V$ be a right $\Gamma$-module. Given a function $\phi: \Delta_0 \rightarrow V$, we say that $\phi$ is \emph{$\Gamma$-invariant} if
\[\phi(D) = \phi|\gamma(D) = \phi(\gamma\cdot D)|\gamma.\]
We define the space of \emph{$V$-valued modular symbols for $\Gamma$} to be the space $\symb_\Gamma(V) \defeq \Hom_{\Gamma}(\Delta_0,V)$ of $\Gamma$-invariant  maps from $\Delta_0$ into $V$.
\end{mdef}
\begin{mrem}
Relating back to the description of modular symbols in the introduction, where we gave motivation for the explicit study of modular symbols in the Bianchi case, we see that using the arguments in \cite{AS86} that there is a canonical isomorphism
\[\h_{\mathrm{c}}^1(\Gamma\backslash\uhs,\mathcal{V}) \cong \symb_\Gamma(V),\]
where $\mathcal{V}$ is the local system corresponding to $V$.
\end{mrem}
Recall that we set $\n$ to be an ideal of $\roi_K$ with $(p)|\n$, and recall the definition of $\Omega_1(\n)$ in Definition \ref{omega}. With our choices of explicit representatives $I_1,...,I_h$ for the class group of $K$ (see the start of Section \ref{BMF}), the discrete subgroups of $\SLt(K)$ corresponding to $\Omega_1(\n)$ are
\begin{align*}\Gamma_i \defeq& \hspace{2pt} \SLt(K) \cap (g_i\Omega_1(\n)\cdot\GLt(\C) g_i^{-1})\\
 =& \left\{\matrd{a}{b}{c}{d} \in \SLt(K): a, d \in 1+ \n, b \in I_i, c \in \n I_i^{-1}\right\}.
 \end{align*}
\begin{mdef}\label{bianchims}\begin{itemize}\item[(i)]
The space of \emph{Bianchi modular symbols of parallel weight $(k,k)$ and level $\Gamma_i$} is defined to be the space
\[\symb_{\Gamma_i}(V_{k,k}^*(\C)) \defeq \Hom_{\Gamma_i}(\Delta_0, V_{k,k}^*(\C))\]
of $V_{k,k}^*(\C)$-valued modular symbols.
\item[(ii)] The space of \emph{Bianchi modular symbols of parallel weight $(k,k)$ and level $\Omega$} is defined to be the space
\[\symb_{\Omega_1(\n)}(V_{k,k}^*(\C)) \defeq \bigoplus_{i=1}^h \symb_{\Gamma_i}(V_{k,k}^*(\C)).\]
\end{itemize}
\end{mdef}

\subsection{Aside on Group Actions}
For the purposes of clarity, we take a slight detour to elaborate on the various group actions we are endowing the polynomial spaces above with. In our previous work, we stated that $V_{2k+2}(\C)$ was an irreducible right-$\SUt(\C)$ module. In the next section, we will need to use a left-action on this space, namely, the one we obtain from the right action above; that is, we define, for $\gamma \in \SUt(\C)$,
\[\gamma \cdot P\left(\binomc{X}{Y}\right) \defeq P|\gamma^{-1}\left(\binomc{X}{Y}\right).\]
In an attempt to avoid confusion, we denote the space by $V_{2k+2}^r(\C)$ when we consider the right action of $\SUt(\C)$, and by $V_{2k+2}^\ell(\C)$ when we consider the corresponding left-action. We also use this notation for the tensor product in the obvious manner. Note that these actions obviously extend to give actions of $\SLt(\C)$ on the spaces involved.

\subsection{Differentials on $\uhs$}
In the sequel, we will require some facts about differentials on $\uhs$; these are listed here for reference.\\
\\
The space $\Omega^1(\uhs,\C)$ of differential 1-forms on $\uhs$ is a 3-dimensional $C^{\infty}(\uhs)$-module, spanned by $dz,dt$ and $d\overline{z}$. The usual action of $\SLt(\C)$ on $\uhs$ induces a left-action of $\SLt(\C)$ on $\Omega^1(\uhs,\C)$ by pull-back. We also have:
\begin{mprop}\label{isodiffspolys}
Let $\Omega^1_0(\uhs,\C)$ be the $\C$-vector space spanned by $dz, dt$ and $d\overline{z}$. There is an isomorphism $\Omega^1_0(\uhs,\C) \rightarrow V_2^\ell(\C)$ of $\SUt(\C)$-modules given by the map sending
\[dz \mapsto A^2,\hsp -dt \mapsto AB, \hsp -d\zbar \mapsto B^2.\]
\end{mprop}
\begin{proof}See \cite{Gha99}, Section 2.2.
\end{proof}
There is a natural map from $\Omega^1(\uhs,\C)$ to $\Omega^1_0(\uhs,\C)$ given by evaluation at $(0,1) \in \uhs$. Combining this map and the isomorphism of Proposition \ref{isodiffspolys} with the action of $\SLt(\C)$ above allows us to define a left-action of $\SLt(\C)$ on $V_2^\ell(\C)$. An explicit check, contained in \cite{Gha99}, shows that this action can be explicitly described as 
\begin{equation}\label{actionofsl2ondiffs}
\gamma\cdot P\binomc{A}{B} = P\left[\frac{1}{|a|^2 + |c|^2}\matrd{\overline{a}}{\overline{c}}{-c}{a}\binomc{A}{B}\right], \hsp \gamma = \matr.
\end{equation}

\subsection{Construction}
Let $\Phi$ be a Bianchi cusp form, giving rise to a collection of $h$ functions $F^i,...,F^h$ on $\GLt(\C)$. To each of these `components' of $\Phi$, we can associate a harmonic $V_{k,k}(\C)$-valued differential $\omega_i$ on $\uhs$, which we can then integrate between cusps to obtain a modular symbol of level $\Gamma_i$. These then combine to give the modular symbol of level $\Omega_1(\n)$ attached to $\Phi$. This process is described in full in \cite{Gha99}, and here we recap the construction. \\
\\
Firstly, note that the Clebsch-Gordan formula gives an injection 
\[\sigma: V_{2k+2}^\ell(\C) \longhookrightarrow V_{k,k}^\ell(\C)\otimes_{\C}V_2^\ell(\C)\]
of $\SUt(\C)$-modules. With this, and using Proposition \ref{isodiffspolys}, we can define a $V_{k,k}^\ell(\C)$-valued differential on $\uhs$ by
\[\omega_{F^i}(g) \defeq g\cdot(\sigma\circ F^i(g)), \hsp g \in \SLt(\C).\]
(It is a simple check to show that $\omega_{F^i}$, viewed as a map from $\SLt(\C)$ into $V_{k,k}^\ell(\C)\otimes_{\C} V_2^\ell(\C)$, is invariant under translations by $\SUt(\C)$ and thus descends to a differential on $\uhs$). We can view this as a $V_{k,k}^r(\C)$-valued differential in the manner described above. This is still not quite what we require; the definition calls for a $V_{k,k}^*(\C)$-valued differential, where here the right action on $V_{k,k}^*(\C)$ is different to that on $V_{k,k}^r(\C)$. This is dealt with via:
\begin{mprop}There is a $\Sigma_0(\roi_L)^2$-equivariant isomorphism
\begin{align*}\eta: V_{k,k}^r(\C) &\longrightarrow V_{k,k}^*(\C)\\
X^qY^{k-q}\Xbar^r\Ybar^{k-r}&\longmapsto \binomc{k}{q}^{-1}\binomc{k}{r}^{-1}\Xc^{k-q}\Yc^q\Xcbar^{k-r}\Ycbar^r.
\end{align*}
\end{mprop}
\begin{proof}This is a simple, but lengthy, check.
\end{proof}

The following is another elementary check, which uses the modularity of $F^i$ and the action of $\SLt(\C)$ on $\Omega^1(\uhs,\C)$:
\begin{mprop}The map
\[\psi_{\f^i} = \psi_{F^i}: \Delta_0 \longrightarrow V_{k,k}^r(\C)\]
defined by
\[\psi_{\f^i}(\{r\}-\{s\}) \defeq \int_r^s \omega_{F^i}\]
is $\Gamma_i$-invariant, and hence $\psi_{\f^i}$ gives a well-defined element of $\symb_{\Gamma_i}(V_{k,k}^r(\C))$. Thus 
\[\phi_{\f^i} \defeq \eta\circ\psi_{\f^i} \in \symb_{\Gamma_i}(V_{k,k}^*(\C)).\]
\end{mprop}
\begin{mdef}
The \emph{modular symbol attached to $\Phi$} is \[\phi_\Phi \defeq (\phi_{\f^1},...,\phi_{\f^h}) \in \symb_{\Omega_1(\n)}(V_{k,k}^*(\C)).\]
\end{mdef}

\subsection{Relation to $L$-values}
\label{explicitcalc}
We now exhibit a link between values of this modular symbol and critical $L$-values of the original Bianchi modular form. To do this, we appeal to \cite{Gha99}, in which he calculates what this symbol looks like more explicitly.
\begin{mprop}\label{cqr}We have, for $a \in K$,
\[\phi_{\f^i}(\{a\} - \{\infty\}) =\sum\limits_{q,r = 0}^k c^i_{q,r}(a) (\Yc-a\Xc)^{k-q}\Xc^q(\Ycbar-\bar{a}\Xcbar)^{k-r}\Xcbar^{r},\]
\[c^i_{q,r}(a) \defeq 2\binomc{2k+2}{k-q+r+1}^{-1}(-1)^{k+r+1}\int_0^\infty t^{q+r}\f^i_{k-q+r+1}(a,t)dt.\]
\end{mprop}
\begin{proof} The bulk of the work is completed in \cite{Gha99}, Section 5.2. We note that as we integrate over vertical paths, we fix $z = a$, and we are not interested in the $dz$ and $d\overline{z}$ terms of the differential.
\end{proof}
This now combines with the integral formula we obtained for the $L$-function $L^i(\Phi,\psi,s)$ in Theorem \ref{integralformula}, where $\psi$ is a Grossencharacter with infinity type 
\[(-u,-v) = \left(\frac{q-r}{2},-\frac{q-r}{2}\right)\]
and conductor $\ff$, to give the required link. We want to set $2s-2 = q+r$, that is, $s = \tfrac{q+r+2}{2}.$ Set $n = k-q+r+1$; then we obtain
\begin{align*}
\sum\limits_{[a]\in\ff^{-1}/\roi_K}\psi(a\ff)^{-1}a^u\overline{a^v}c_{q,r}^i(a) = (-1)^{k+r+1}2\smallbinomc{2k+2}{n}^{-1}A(i,n,\psi,\tfrac{q+r+2}{2})^{-1}L^i(\Phi,\psi,\tfrac{q+r+2}{2}).
\end{align*}
Some cancellation and the explicit form for $A(i,n,\psi,\tfrac{q+r+2}{2})$ now gives
\begin{align*}
\sum\limits_{\substack{[a]\in\ff^{-1}/\roi_K\\ ((a)\ff, \ff) = 1}}&\psi(a\ff)^{-1}a^u\overline{a^v}c_{q,r}^i(a)\\
 &= \left[\psi(t_i)|t_i|_f^{\frac{q+r}{2}}\cdot\frac{(-1)^{q+r+2+k}2(2\pi i)^{q+r+2}}{|\delta|^{q+r+2}\itGamma(q+1)\itGamma(r+1)w\tau(\psi^{-1})}\right]^{-1}L^i(\Phi,\psi,\tfrac{q+r+2}{2}).
\end{align*}
This gives us a link between modular symbols and $L$-values. Indeed, combining the results above for each $i$, we have
\begin{mprop}In the set-up of above, we have
\begin{align*}
L(\Phi,&\psi,\tfrac{q+r+2}{2}) =\\
 &\left[\frac{(-1)^{k}2(-2\pi i)^{q+r+2}}{|\delta|^{q+r+2}q!r!w\tau(\psi^{-1})}\right]\sum\limits_{i=1}^h\bigg[\psi(t_i)|t_i|_f^{\frac{q+r}{2}}\sum\limits_{\substack{[a]\in\ff^{-1}/\roi_K\\ ((a)\ff, \ff) = 1}}\psi(a\ff)^{-1}a^u\overline{a^v}c_{q+r}^i(a)\bigg].
\end{align*}
\end{mprop}

\subsection{$L$-Functions as Functions on Characters}\label{lfnchar}
It is convenient to rephrase the above; instead of seeing the $L$-function as a complex function of one variable, we think of it as a function on Grossencharacters by setting 
\[L(\Phi,\psi) = L(\Phi,\psi,1).\]

Let $\psi' \defeq \psi|\cdot|^{\tfrac{q+r}{2}}.$ This is now a Grossencharacter of conductor $\ff$ and infinity type $(q,r)$. To see that there is a relation of Gauss sums between $\psi$ and $\psi'$, we need some new notation. To any ideal $I$ of $K$, we associate an idele $x_I$ in the following way: if $(\alpha)$ is principal, then define
\[(x_{(\alpha)})_v \defeq \left\{\begin{array}{ll}\alpha &: v = \pri\text{ is finite, and } \pri|(\alpha),\\
1 &: \mathrm{otherwise}.\end{array}\right. .\]
Then for $I = \alpha I_i,$ define $x_I \defeq t_ix_{(\alpha)}.$ Note that if $I$ is coprime to $\ff$, then $\psi(I) = \psi(x_I).$ Now, we see that
\[\tau((\psi')^{-1}) = |\delta|^{q+r}|x_{\ff}|_f^{-\frac{q+r}{2}}\tau(\psi^{-1}),\]
which is proved using the easily seen identity $|x_{(a)\ff}|_{\A_K} = |x_{(a)\ff}|_f= |a|_f|x_{\ff}|_f = |a|_{\infty}^{-2}|x_{\ff}|_f$. Accordingly, we see that
\begin{align*}L(\Phi,\psi')& = \left[\frac{(-1)^{k}2(-2\pi i)^{q+r+2}}{|\delta|^{2}|x_{\ff}|_f^{\frac{q+r}{2}}q!r!w\tau((\psi')^{-1})}\right]
\sum\limits_{i=1}^h\bigg[\psi'(t_i)\sum\limits_{\substack{[a]\in\ff^{-1}/\roi_K\\ ((a)\ff, \ff) = 1}}\psi'(a\ff)^{-1}|x_{(a)\ff}|^{\frac{q+r}{2}}a^u\overline{a^v}c_{q,r}^i(a)\bigg]\\
&= \left[\frac{(-1)^{k}2(-2\pi i)^{q+r+2}}{Dq!r!w\tau((\psi')^{-1})}\right]
\sum\limits_{i=1}^h\bigg[\psi'(t_i)\sum\limits_{\substack{[a]\in\ff^{-1}/\roi_K\\ ((a)\ff, \ff) = 1}}\psi'(a\ff)^{-1}(\psi'_\infty)^{-1}(a)c_{q,r}^i(a)\bigg],
\end{align*}
where here we have used that $|\delta|^2 = D$. This simplifies further; indeed, for any Grossencharacter $\psi$, an explicit check shows that when $a \in \ff^{-1}/\roi_K$ with $(a\ff,\ff) = 1$, we have
\[\psi(a\ff)^{-1}\psi_\infty(a)^{-1} = \psi(x_{\ff})^{-1}\psi_{\ff}(ax_{\ff}),\]
so that this becomes
\begin{align*}L(\Phi,\psi') = \left[\frac{(-1)^{k}2(-2\pi i)^{q+r+2}}{\psi'(x_{\ff})Dq!r!w\tau((\psi')^{-1})}\right]
\sum\limits_{i=1}^h\bigg[\psi'(t_i)\sum\limits_{\substack{[a]\in\ff^{-1}/\roi_K\\ ((a)\ff, \ff) = 1}}\psi'_{\ff}(ax_{\ff})c_{q,r}^i(a)\bigg],
\end{align*}
where $\psi'$ is a Grossencharacter of conductor $\ff$ and infinity type $(q,r)$ with $0\leq q,r \leq k$.\\
\\
We make one further change with the aim of massaging this formula into something a little nicer; namely, we renormalise, using the Deligne $\Gamma$-factor at infinity. Define
\[\Lambda(\Phi,\psi,t) \defeq \frac{\itGamma(q+t)\itGamma(r+t)}{(2\pi i)^{q+t}(2\pi i)^{r+t}}L(\Phi,\psi,t),\]
from which we deduce that:
\begin{mthm}\label{lfunctionmodsymb}Let $K/\Q$ be an imaginary quadratic field of discriminant $D$ and $w$ be the size of $\roi_K^\times$, and let $\Phi$ be a cuspidal Bianchi modular form of parallel weight $(k,k)$ and level $\Omega_1(\n)$ with renormalised $L$-function $\Lambda(\Phi,\psi)$. For a Grossencharacter $\psi$ of $K$ of conductor $\ff$ and infinity type $0 \leq(q,r)\leq (k,k)$, we have
\begin{equation}\label{finallformula}
\Lambda(\Phi,\psi) = \left[\frac{(-1)^{k+q+r}2\psi_{\ff}(x_{\ff})}{\psi'(x_{\ff})Dw\tau((\psi')^{-1})}\right]\sum\limits_{i=1}^h\bigg[\psi(t_i)\sum\limits_{\substack{[a]\in\ff^{-1}/\roi_K\\ ((a)\ff, \ff) = 1}}\psi_{\ff}(a)c_{q,r}^i(a)\bigg],
\end{equation}
where $c_{q,r}^i(a)$ is as defined in Proposition \ref{cqr}.
\end{mthm}

\section{Overconvergent Bianchi Modular Symbols}
\label{ombssec}
Overconvergent modular symbols are modular symbols that take values in a space of $p$-adic distributions. In a similar style to the rest of the paper, we begin this section by defining these spaces in their most natural way - in the process, directly generalising the work of Pollack and Stevens - before passing to a more workable and explicit description. We'll then write down a filtration on the space of overconvergent Bianchi modular symbols of fixed weight and level, putting us in the situation of the next section and allowing us to prove a control theorem for these objects.\\
\\
In the following, let $F$ be either $\Q$ or an imaginary quadratic field with ring of integers $\roi_F$. We can write down the polynomial space associated to modular symbols (either $V_k(\Cp)$ over $\Q$ or $V_{k,k}(\Cp)$ over an imaginary quadratic field) in a more intrinsic way, namely as the space $V_k(\roi_F\otimes_{\Z}\Zp)$ of polynomial functions on $\roi_F\otimes_{\Z}\Zp$ of degree less than or equal to $k$ in each variable. The space $V_k(\roi_F\otimes_{\Z}\Zp)$ is a subspace of the space $\A(\roi_F\otimes_{\Z}\Zp)$ of rigid analytic functions on $\roi_F\otimes_{\Z}\Zp$. We then find that, by dualising the inclusion $V_k(\roi_F\otimes_{\Z}\Zp) \hookrightarrow \A(\roi_F\otimes_{\Z}\Zp)$, we have a surjection from the space $\D(\roi_F\otimes_{\Z}\Zp)$ of rigid analytic distributions on $\roi_F\otimes_{\Z}\Zp$ to the space $V_k^*(\roi_F\otimes_{\Z}\Zp)$, giving a surjective specialisation map. To see that we can define modular symbols with values in $\D(\roi_F\otimes_{\Z}\Zp)$, we consider the semigroup 
\[\Sigma_0(\roi_F\otimes_{\Z}\Zp) \defeq \left\{\matr \in M_2(\roi_F\otimes_{\Z}\Zp): c \in p\roi_F\otimes_{\Z}\Zp, a \in (\roi_F\otimes_{\Z}\Zp)^\times, ad-bc \neq 0\right\}.\]
(Note that when $p$ ramifies as $\pri^2$ in $F$, we can consider the larger group coming from the condition $c \in \pi_{\pri}\roi_F\otimes_{\Z}\Zp$, where $\pi_{\pri}$ is a uniformiser at $\pri$. We make no further mention of this here, however). This has a natural action on $\A(\roi_F\otimes_{\Z}\Zp)$, depending on $k$, defined as follows: for $\gamma = \smallmatrd{a}{b}{c}{d} \in \Sigma_0(\roi_F\otimes_{\Z}\Zp)$ and $f \in \A(\roi_F\otimes_{\Z}\Zp)$, define
\[\gamma \cdot_k f(x) = (a+cx)^kf\left(\frac{b+dx}{a+cx}\right).\]
Using this we obtain an action of suitable discrete subgroups of $\SLt(F)$ on the distribution space, allowing us to use it as a value space. We also get a Hecke action on the resulting modular symbols. All of this directly generalises the work of Pollack and Stevens in the case $F = \Q$, and with suitable small adjustments would generalise further to the case of arbitrary number fields. \\
\\
In our setting, we can describe these distribution spaces more explicitly. Later, when $p$ is split or inert, we will show that the space $\A(\roi_F\otimes_{\Z}\Zp)$ can be thought of as just being the space of rigid analytic functions on $\Zp^2$, or, in the weak topology, the completed tensor product of two copies of the space of rigid analytic functions on one copy of $\Zp$. The group $\Sigma_0(\roi_F\otimes_{\Z}\Zp)$ naturally injects into two copies of the semigroup $\Sigma_0(p)$ as written down by Pollack and Stevens (we will call this semigroup $\Sigma_0(\roi_L)$). When $p$ is not split, we are carrying around some redundant information with this approach; in particular, the useful information given by the action of $\Sigma_0(p)^2$ is entirely determined by the action of one of the components. There are significant advantages to using this more explicit approach, however. The spaces have nice descriptions that are easy to work with and allow us to generalise the filtration proof of Greenberg (\cite{Gre07}) to the imaginary quadratic case. In the remainder of this section, we work with this explicit approach exclusively. 
\subsection{Notation and Preliminaries}
\label{notation}
\begin{mnotnum}\label{not}
Throughout, as before, $K/\Q$ denotes an imaginary quadratic field. Let $p$ be a rational prime with $p\roi_K = \prod \pri^{e_{\pri}}$ and define $f_{\pri}$ to be the residue class degree of $\pri$. Note that $\sum e_{\pri}f_{\pri} = 2$. Fix an embedding $\overline{\Q} \hookrightarrow \overline{\Qp}$; then for each prime $\pri|p$, we have $e_{\pri}f_{\pri}$ embeddings $K_{\pri} \hookrightarrow \overline{\Qp}$, and combining these for each prime, we get an embedding 
\[\sigma : K\otimes_{\Q}\Qp \longhookrightarrow \overline{\Qp}\times\overline{\Qp}.\]
If $a \in K\otimes_{\Q}\Qp$, we write $\sigma_1(a)$ and $\sigma_2(a)$ for the projection of $\sigma(a)$ onto the first and second factors respectively. Now, let $L$ be a finite extension of $\Qp$ such that the image of $\sigma$ lies in $L^2$. We equip $L$ with a valuation $v$, normalised so that $v(p) = 1$, and denote the ring of integers in $L$ by $\roi_L$, with uniformiser $\pi_L$. Then in fact, for each integral ideal $I$ of $K$ coprime to $(p)$, we have 
\begin{align}\label{embedding}\sigma:I^{-1} \longhookrightarrow \roi_K\otimes_{\Z}\Zp \longhookrightarrow \roi_L \times \roi_L.\end{align}
In the obvious way, we then have an embedding
\[\left\{\matrd{a}{b}{c}{d}: a,b,c,d \in I^{-1}, ad-bc \neq 0\right\} \longhookrightarrow \GLt(\roi_L)\times\GLt(\roi_L).\]

\end{mnotnum}

Note that in particular this means that the groups in the following definition - and hence the groups $\Gamma_i = \Gamma_{I_i}$ of previous sections - embed into the right hand side. Whilst the groups $\Gamma_I$ are implicitly dependent on the ideal $\n$, since we work exclusively at level $\n$ throughout this paper, for convenience we drop $\n$ from the notation.
\begin{mdef}\label{gammaidef}
Let $\n$ be an ideal of $\roi_K$ with $(p)|\n$, and let $I$ be an integral ideal of $K$ that is coprime to $\n$. Define the \emph{twist of $\Gamma_1(\n)$ by $I$} to be
\[\Gamma_I \defeq \left\{\matrd{a}{b}{c}{d} \in \GLt(K): a,d \in 1 + \n, b \in I, c \in \n I^{-1}\right\}.\]
\end{mdef}
The embedding above will be used in the sequel to define the action of each $\Gamma_i$ as well a Hecke action on suitable modular symbol spaces. We'll define some monoid $\Sigma_0(\roi_L) \leq \GLt(\roi_L)$ and an action of $\Sigma_0(\roi_L)\times\Sigma_0(\roi_L)$. Every matrix whose action we study will have image in $\Sigma_0(\roi_L)^2$ under the embedding above. Thus in proving facts about the action of $\Sigma_0(\roi_L)\times\Sigma_0(\roi_L),$ we'll encapsulate everything we'll later need \emph{regardless} of the splitting behaviour of $p$ in $\roi_K$.
\begin{mrem}For our purposes, we may need to take $L$ to be larger than this. Let $\Phi$ be a cuspidal Bianchi eigenform with Fourier coefficients $c(I,\Phi)$, normalised so that $c(\roi_K,\Phi) = 1$. Then the Fourier coefficients are algebraic. In particular, when studying the action of the $U_p$ operator we may need to consider eigenvalues living in the number field $F \defeq K(\{c(\pri,\Phi): \pri|p\})$. As we can easily enlarge $L$ to contain all possible embeddings of the completions of this field into $\overline{\Q_p}$, we will henceforth assume that these eigenvalues live in $L$.
\end{mrem}
\subsection{Overconvergent Modular Symbols}
To define overconvergent modular symbols, Stevens used spaces from $p$-adic analysis. For more details on the results here, including $p$-adic function and distribution spaces as well as the completed tensor product, see \cite{Colm10}.
\begin{mdef}[Modules of Values for Overconvergent Modular Symbols]
Let $R$ be either a $p$-adic field or the ring of integers in a finite extension of $\Qp$, and let $\A(R)$ be the ring of \emph{rigid analytic functions on the closed unit disc defined over $R$,} that is, the ring
\[\A(R) = \left\{\sum_{n\geq0} a_nx^n: a_n \in R, a_n \rightarrow 0 \text{ as $n$ tends to $\infty$}\right\}.\]
For $R = L$, this is a $L$-Banach space with the sup norm. We write $\A_2(R)$ for the completed tensor product $\A(R)\ctp_R\A(R).$ We let $\D(R) = \Hom_\text{cts}(\A(R),R)$, the space of \emph{rigid analytic distributions over $R$}, be the topological dual of $\A(R)$, and analogously we let $\D_2(R) = \Hom_\text{cts}(\A_2(R),R)$ be the topological dual of $\A_2(R)$.
\end{mdef}
We have the following identification, as promised in the introduction to this section.
\begin{mprop}\label{tensorisom}If $p$ is split or inert in $K$, then the natural inclusion $\A_2(L) \hookrightarrow \A(\roi_K\otimes_{\Z}\Zp)$, given by restriction, is an isomorphism.
\end{mprop}
\begin{proof}Suppose $p$ is split. Then, up to an ordering of the primes above $p$, there is a natural identification of $\roi_K\otimes_{\Z}\Zp$ with $\Zp^2$, and the result follows.\\
\\
Suppose instead that $p$ is inert. Then a rigid analytic function $f$ on $\roi_K\otimes_{\Z}\Zp$ takes the form
\[f: x \longmapsto \sum_{r,s}a_{r,s}x^r\overline{x}^s, \hspace{12pt} a_{r,s} \in L,\]
where the coefficients lie in $L$ as we took $L$ large enough to contain $K\otimes_{\Q}\Qp$. Moreover, rigid analyticity means that the coefficients $a_{r,s}$ tend to zero as $r+s$ tends to infinity. Accordingly, we can identify $f$ with an element $\tilde{f}(X,Y) = \sum_{r,s}a_{r,s}X^rY^s$ of the Tate algebra $L\langle X,Y\rangle$ in two variables defined over $L$. But this can be viewed as the space of rigid analytic functions on $\Zp^2$ defined over $L$, and the result follows.
\end{proof}
\begin{mrems}\begin{itemize}
\item[(i)]When $p$ is inert, the fact that we can make this identification relies on our functions being rigid analytic. If instead we considered functions that were locally analytic of order 1, such an identification does not hold. This is because any rigid analytic function on $\roi_K\otimes_{\Z}\Zp$ can be extended to $\roi_{\Cp}$, whereas the same is \emph{not} true of locally analytic functions.
\item[(ii)] Furthermore, this identification does not carry over immediately to the case $p$ ramified. Indeed, consider the function $x \mapsto \sum_{r,s}\pi_{\pri}^n([x-\overline{x}]/\pi_{\pri})^{n^2}$. This is a rigid analytic function on $\roi_K\otimes_{\Z}{\Zp}$ that cannot be written as an element of the Tate algebra over $L$. However, we will work with the more explicit definition in the case $p$ ramified as well, as ultimately the construction still gives us a $p$-adic $L$-function in this case. In this case, the relevant distribution space for constructing the $p$-adic $L$-function, $\D(\roi_K\otimes_{\Z}\Zp)$, is naturally a subspace of $\D_2(L)$; indeed, in Proposition \ref{locanalyticprop}, we show that the overconvergent modular symbol we construct lives in this smaller space (and, indeed, in the even smaller space of \emph{locally analytic distributions}).
\end{itemize}
\end{mrems}
Having defined these spaces of distributions, our primary spaces of interest, we immediately give two alternate descriptions that are easier to work with.
\begin{mdef}
\begin{itemize}
\item[(i)] Let $\mu \in \D$ be a distribution. Define the \emph{moments} of $\mu$ to be the values $(\mu(x^i))_{i\geq 0},$ noting that these values totally determine the distribution since the span of the $x^i$ is dense in $\A$.
\item[(ii)] Let $\mu \in \D_2$ be a two variable distribution. Define the \emph{moments} of $\mu$ analogously to be the values $(\mu(x^iy^j))_{i,j \geq 0}$.
\end{itemize}
\end{mdef}
The following simple proposition gives a simple description of the distribution spaces.
\begin{mprop}
By taking the moments of a distribution, we can identify $\D(L)$ with the set of bounded sequences in $L$, and $\D_2(L)$ with the set of doubly indexed bounded sequences in $L$.
\end{mprop}

\begin{mremsnum}\label{tensor}
This identification with bounded sequences means that we have $\D_2(L) \cong \D_2(\roi_L) \otimes_{\roi_L} L,$ where $\D_2(\roi_L)$ is the subspace of $\D_2(L)$ consisting of distributions with integral moments.
\end{mremsnum}

\subsection{The Action of $\Sigma_0(\roi_L)^2$ and Hecke Operators}
\label{actionofsig}
We want to equip our spaces with an action of suitable congruence subgroups and a Hecke action. Recall the definition of $\Sigma_0(\roi_K\otimes_{\Z}\Zp)$ in the introduction to this section, and its left weight $k$ action 
\[(\gamma \cdot_k f)(x) = (cx + a)^kf\left(\frac{dx+b}{cx + a}\right), \hsp \gamma  = \matr\]
on rigid analytic functions on $\roi_K\otimes_{\Z}\Zp$. To obtain a description of this action in the more explicit setting, define 
\[\Sigma_0(\roi_L) \defeq \left\{\matr \in M_2(\roi_L): p\mid c, \hspace{3pt}(a,p) = 1, ad-bc \neq 0\right\}.\]
Note here that $\Sigma_0(\roi_K\otimes_{\Z}\Zp)$ embeds inside $\Sigma_0(\roi_L)^2$ in a natural way; indeed, if $p$ is split, then $\roi_K\otimes_{\Z}\Zp$ is isomorphic to $\Zp^2$ and the inclusion is clear. If $p$ is inert or ramified then the natural embedding $\sigma: \roi_K\otimes_{\Z}\Zp \hookrightarrow \roi_L^2$ of equation (\ref{embedding}) extends in the natural way to an embedding of the matrix groups. Now suppose $R$ contains $\roi_L$. Then the left weight $k$ action of $\Sigma_0(\roi_K\otimes_{\Z}\Zp)$ on $\A(\roi_K\otimes_{\Z}\Zp)$ extends in the obvious way to give a weight $(k,k)$ and, more generally, a weight $(k,\ell)$ action of $\Sigma_0(\roi_L)^2$ on $\A_2(R)$ by
\[((\gamma_1,\gamma_2)\cdot_{(k,\ell)}f)(x,y) = (c_1x+a_1)^k(c_2y+a_2)^\ell f\left(\frac{d_1x+b_1}{c_1x+a_1},\frac{d_2y+b_2}{c_2y+a_2}\right),\]
where $\gamma_i = \smallmatrd{a_i}{b_i}{c_i}{d_i}.$ This gives rise to a right weight $(k,\ell)$ action of $\Sigma_0(\roi_L)^2$ on $\D_2(R)$ defined by
\[\mu|_{k,\ell}(\gamma_1,\gamma_2)(f) = \mu((\gamma_1,\gamma_2)\cdot_{(k,\ell)}f).\]
When talking about these spaces equipped with the weight $(k,\ell)$ action, we denote them by $\A_{k,\ell}(R)$ and $\D_{k,\ell}(R)$ respectively.
\begin{mremsnum}
\begin{itemize}
\item[(i)] Importantly, for $p$ split or inert, the isomorphism of Proposition \ref{tensorisom} respects the actions of $\Sigma_0(\roi_K\otimes_{\Z}\Zp)$ on the left hand side and $\Sigma_0(\roi_L)^2$ on the right hand side.
\item[(ii)]\label{actiononvkl}Note that the subspace $V_{k,\ell}(R)$ of $\A_2(R)$ is stable under the action of $\Sigma_0(\roi_L)$, and hence it inherits a left action of $\Sigma_0(\roi_L)$. This is the action we earlier defined in Definition \ref{actionofsltc}.
\end{itemize}
\end{mremsnum}

Note that for any ideal $I$ coprime to $\n$, we have a right action of $\Gamma_I$ (as in Definition \ref{gammaidef}) on the space $\D_{k,\ell}(L).$
\begin{mdef}\begin{itemize}
\item[(i)] Define the space of \emph{overconvergent Bianchi modular symbols for $K$ of weights ($k,\ell)$ and level $\Gamma_I$ with coefficients in $L$} to be
\[\symb_{\Gamma_I}(\D_{k,\ell}(L)) \defeq \Hom_{\Gamma_I}(\Delta_0, \D_{k,\ell}(L)).\]
\item[(ii)] Recall the definitions of $\Omega_1(\n)$ from Definition \ref{omega} and $\Gamma_i$ in equation (\ref{gammai}). Define the space of \emph{overconvergent Bianchi modular symbols for $K$ of weights ($k,\ell)$ and level $\Omega_1(\n)$ with coefficients in $L$} to be
\[\symb_{\Omega_1(\n)}(\D_{k,\ell}(L)) \defeq \bigoplus_{i=1}^h\symb_{\Gamma_i}(\D_{k,\ell}(L)).\]
\end{itemize}
\end{mdef}
For a right $\Sigma_0(\roi_L)^2$-module $V$, we define a Hecke action on the adelic space $\symb_{\Omega_1(\n)}(V) \defeq \bigoplus_{i=1}^h\symb_{\Gamma_i}(V)$ as follows. For a prime $\pri$ dividing $(p)$, for each $i \in \{1,...,h\}$ there is a unique $j_i \in \{1,...,h\}$ such that
\[\pri I_i = (\alpha_i)I_{j_i},\]
for $\alpha_i \in K$. Then the $U_{\pri}$ operator is
\[(\phi_1,...,\phi_h)|U_{\pri} \defeq \left(\phi_{j_1}\bigg|\left[\Gamma_{j_1}\matrd{1}{0}{0}{\alpha_1}\Gamma_1\right],...,\phi_{j_h}\bigg|\left[\Gamma_{j_h}\matrd{1}{0}{0}{\alpha_h}\Gamma_h\right]\right).\]
We can work out the double coset operators explicitly to be given by
\[\left[\Gamma_{j_i}\matrd{1}{0}{0}{\alpha_i}\Gamma_i\right] = \sum_{a\newmod{\pri}}\matrd{1}{a}{0}{\alpha_i}\]
using the usual methods. Note that if $n$ is an integer such that $\pri^n = (\sigma)$ is principal, then this action becomes significantly simpler; namely, we just act on each component, with no permuting, via the sum $\sum_{a\newmod{\pri^n}}\smallmatrd{1}{a}{0}{\sigma}$. Because of this much simpler description, in the sequel we much prefer to use a principal power of $\pri$ and prove results using just one component at a time. \\
\\
Ideally, we'd prefer to work with integral distributions.

\begin{mlem}\label{finitegen}Let $K$ be a number field and $I$ be an ideal coprime to $\n$. Then $\mathrm{Div}^0(\Proj(K))$ is a finitely generated $\Z[\Gamma_I]$-module.
\end{mlem}
\begin{proof}
This follows from the fact that $\Gamma_I$ is a finitely generated group (see, for example, \cite{Swa71}) and that it has finitely many cusps.
\end{proof}

\begin{mprop}\label{tensoring}Let $\Gamma_I$ be as above, and let $D$ have the structure of both a $\roi_L$-module and a right $\Gamma_I$-module. Then we have
\[\symb_{\Gamma_I}(D \otimes_{\roi_L} L) \cong \symb_{\Gamma_I}(D) \otimes_{\roi_L} L.\]
\end{mprop}
\begin{proof}
Let $\phi \in \symb_{\Gamma_I}(D \otimes_{\roi_L} L).$ Using Lemma \ref{finitegen}, take a finite set of generators $\alpha_1, ..., \alpha_n$ for $\Delta_0$ as a $\Z[\Gamma_I]$-module. We can find some element $c \in \roi_L$ such that $c\phi(\alpha_i) \in D$ for each $i$. But then it follows immediately that $c\phi \in \symb_{\Gamma_I}(D)$, and the result follows.
\end{proof}
With this structure in place, we can now work with the space $\D_{k,\ell}(\roi_L)$.

\subsection{Finite Approximation Modules}
\begin{mrem}From now on, we work with parallel weights, i.e. we consider only $k = \ell$ and use the space $\D_{k,k}(R).$ There are no classical cuspidal Bianchi modular forms at non-parallel weights, so in proving a control theorem in the spirit of Stevens' work, it suffices to exclude the case $k \neq \ell$. We'll also focus on looking at the space of Bianchi modular symbols one component at a time. Henceforth, to this end, $\Gamma$ will denote one of the $\Gamma_i$ for $i \in \{1,...,h\}$.
\end{mrem}
In the one variable case, in \cite{Gre07} Matthew Greenberg gave an alternative proof of Stevens' control theorem using \emph{finite approximation modules}, defining a $\Sigma_0(\roi_L)$-stable filtration of $\D_k(L)$ and then lifting modular symbols through this filtration. We aim to mimic this. First we recap Greenberg's filtration, recasting it slightly to make it more favourable for our generalisation. He defines:
\begin{mdef}\label{onevarfilt}
\begin{itemize}
\item[(i)]$\mathcal{F}^N\D_k(\roi_L) \defeq \{\mu \in \D_k(\roi_L): \mu(x^i) \in \pi_L^{N-i}\roi_L\},$ and
\item[(ii)] $\D_k^0(\roi_L) = \{\mu\in D_k(\roi_L): \mu(x^i) = 0$ for $0\leq i \leq k\}$. Note that this is the kernel of the natural map $\D_k(\roi_L) \rightarrow V_k^*(\roi_L)$ obtained by dualising the inclusion $V_k(\roi_L) \hookrightarrow \A_k(\roi_L)$. 
\end{itemize}
\end{mdef}
To define our own filtration, we take a similar route, imposing suitable conditions on the moments of distributions.
\begin{mdef}
Define:
\begin{itemize}
\item[(i)] $\mathcal{F}^N\D_{k,k}(\roi_L) \defeq \{\mu \in \D_{k,k}(\roi_L): \mu(x^iy^j) \in \pi_L^{N -i - j}\roi_L\}.$
\item[(ii)] $\D_{k,k}^0(\roi_L) \defeq \{\mu \in \D_{k,k}(\roi_L): \mu(x^iy^j) = 0$ for $0\leq i,j \leq k\}.$
\item[(iii)] $F^N\D_{k,k}(\roi_L) \defeq \mathcal{F}^N\D_{k,k}(\roi_L) \cap \D_{k,k}^0(\roi_L).$
\end{itemize}
\end{mdef}

\begin{mprop}\label{sigstable}
This filtration is $\Sigma_0(\roi_L)^2$-stable.
\end{mprop}
\begin{proof}
There is an obvious switching map $s: \D_{k,\ell}(\roi_L) \rightarrow \D_{\ell,k}(\roi_L)$. Thus it suffices to prove the result for elements of form $(\gamma, I_2)$ for $\gamma \in \Sigma_0(\roi_L)$, as the action of a more general element  can be described as
\[\mu\bigg|_k(\gamma_1,\gamma_2) = s^{-1}\left[s\left(\mu\bigg|_k(\gamma_1,I_2)\right)\bigg|_k(\gamma_2,I_2)\right].\]
To each two-variable distribution $\mu \in \D_{k,k}(\roi_L),$ associate a family of distributions $\{\mu_j \in \D_k(\roi_L)\}$ by defining the moments of $\mu_j$ to be
\[\mu_j(x^i) = \mu(x^iy^j).\]
Then note that we have
\[\mu \in F^N\D_{k,k}(\roi_L) \iff \mu_j \in \left\{\begin{array}{ll}\mathcal{F}^{N-j}\D_k(\roi_L) \cap \D_k^0(\roi_L) & : 0\leq j \leq k \\ \mathcal{F}^{N-j}\D_k(\roi_L) & : j > k\end{array}\right.,\]
where the condition must hold for all $j \geq 0$. The result we require then follows from the observation that 
\[\mu\bigg|_k(\gamma, I_2)(x^iy^j) \hspace{3pt}=\hspace{3pt} \mu_j\bigg|_k\gamma(x^i)\]
combined with the stability (in the one variable case) of each of the modules $\D_k^0(\roi_L)$ and $\mathcal{F}^{N-j}\D_k(\roi_L)$  under the action of $\Sigma_0(\roi_L)$.
\end{proof}
In particular, this result shows that we have a collection of $\Sigma_0(\roi_L)^2$-modules 
\[A^N\D_{k,k}(\roi_L)\defeq \frac{\D_{k,k}(\roi_L)}{F^N\D_{k,k}(\roi_L)},\]
with action inherited from $\D_{k,k}(\roi_L)$, and where this is well-defined since the $F^N\D_{k,k}(\roi_L)$ are $\Sigma_0(\roi_L)^2$-stable. Furthermore, we see that
\begin{align*}A^N\D_{k,k}(\roi_L) &\cong \roi_L^{(k+1)^2} \times T\\
&\cong V_{k,k}^*(\roi_L) \times T,
\end{align*}
where $T$ is some finite product of copies of $\Z/p, \Z/p^2,$ and so on up to  $\Z/p^{N-k-1}$. In particular, we also have $A^0\D_{k,k}(\roi_L) \cong \cdots \cong A^{k+1}\D_{k,k}(\roi_L) \cong V_{k,k}^*(\roi_L)$. We also have $\Sigma_0(\roi_L)^2$-equivariant projection maps $\pi^N$ from $\D_{k,k}(\roi_L)$ to $A^N\D_{k,k}(\roi_L)$, and we see that $\pi_0$ is in fact the  map
\[\pi^0: \D_{k,k}(\roi_L) \longrightarrow V_{k,k}^*(\roi_L)\]
that gives rise to the ($\Sigma_0(\roi_L)^2$-equivariant) specialisation map
\[\rho^0: \symb_\Gamma(\D_{k,k}(\roi_L)) \longrightarrow \symb_\Gamma(V_{k,k}^*(\roi_L)).\]
Further to this, for each $M \geq N$, we have a $\Sigma_0(\roi_L)^2$-equivariant map 
\[\pi^{M,N}: A^M\D_{k,k}(\roi_L) \rightarrow A^N\D_{k,k}(\roi_L)\]
 given by projection (and hence also maps $\rho^{M,N}$). The projection maps are all compatible in the obvious ways. Thus we get an inverse system of modules, and it's straightforward to see that
\begin{align}\label{inverselimit}
\D_{k,k}(\roi_L) \cong \lim\limits_{\longleftarrow}A^N\D_{k,k}(\roi_L).
\end{align}
It follows from equation (\ref{inverselimit}) that 
\begin{align}\label{inverselimit2}\symb_\Gamma(\D_{k,k}(\roi_L)) \cong \lim\limits_{\longleftarrow}\symb_\Gamma(A^N\D_{k,k}(\roi_L)).\end{align}

\subsection{Moments of Functionals on Polynomials}
Given an element $f \in V_{k,k}^*(\roi_L)$, define the \emph{moments of $f$} to be the quantities $f(x^my^n)$ for $0 \leq m,n \leq k$, analogously to the moments of more general distributions from previously. Any such $f$ is entirely determined by its moments.
\begin{mdef}Let $\lambda \in L^*$. Define
\[V_{k,k}^\lambda(\roi_L) \defeq \left\{f\in V_{k,k}^*(\roi_L): f(x^iy^j) \in \lambda p^{-(i+j)}\roi_L, 0 \leq i+j \leq \lfloor v(\lambda)\rfloor\right\}.\]
\end{mdef}

\begin{mprop}\label{vlstable}
$V_{k,k}^\lambda(\roi_L)$ is $\Sigma_0(\roi_L)^2$-stable.
\end{mprop}
\begin{proof}Argue as before, or see \cite{Gre07}.
\end{proof}

The next result is a technical lemma describing the action of certain matrices in $\Sigma_0(\roi_L)^2$. It gives us nice properties of the $U_p$ operator.
\begin{mlem}\label{tech}\begin{itemize}
\item[(i)] Let $\mu \in \D_{k,k}(\roi_L)$ be such that $\pi^0(\mu) \in V_{k,k}^\lambda(\roi_L).$ Then, for $a_i \in \roi_L$, we have
\[\mu\bigg|_k\left[\matrd{1}{a_1}{0}{p},\matrd{1}{a_2}{0}{p}\right] \in \lambda\D_{k,k}(\roi_L).\]
\item[(ii)] Let $\mu \in F^N\D_{k,k}(\roi_L),$ and suppose $v(\lambda) < k+1.$ Then
\[\mu\bigg|_k\left[\matrd{1}{a_1}{0}{p},\matrd{1}{a_2}{0}{p}\right] \in \lambda F^{N+1}\D_{k,k}(\roi_L).\]
\end{itemize}
\end{mlem}
\begin{proof}Take some $\mu \in \D_{k,k}(\roi_L).$ Then
\[\mu\bigg|_k\left(\matrd{1}{a_1}{0}{p},\matrd{1}{a_2}{0}{p}\right)(x^my^n) = \mu((a_1+px)^m(a_2+py)^n)\]
\[=\sum\limits_{i=0}^m\sum\limits_{j=0}^n\binomc{m}{i}\binomc{n}{j}a_1^{m-i}a_2^{n-j}p^{i+j}\mu(x^iy^j).\]
\begin{itemize}
\item[(i)] Suppose that $\pi^0(\mu)$ lies in $V_{k,k}^\lambda(\roi_L),$ so that $\mu(x^iy^j) \in \lambda p^{-(i+j)}\roi_L$ for any $i + j \leq\lfloor v(\lambda)\rfloor$. It follows that each term of the sum above lies in $\lambda\roi_L$, and hence the result follows. If instead $i+j > v(\lambda)$, then the result follows as $\mu(x^iy^j) \in \roi_L$.
\item[(ii)] Now suppose $\mu \in F^N\D_{k,k}(\roi_L).$ Again considering the sum above, the terms where $i,j \leq k$ vanish. If $i+j > k$, then 
\[i+j \geq k+1 > v(\lambda),\]
since $\lambda$ has $p$-adic valuation $<k+1$. As $p^{i+j}$ and $\lambda$ are divisible by integral powers of $\pi_L$, it follows that $p^{i+j} \in \pi_L\lambda\roi_L.$ Hence, as $\mu(x^iy^j) \in \pi_L^{N-i-j}\roi_L$, it follows that
\[p^{i+j}\mu(x^iy^j) \in \lambda \pi_L^{(N+1)-i-j}\roi_L,\]
which completes the proof.\qedhere
\end{itemize}
\end{proof}

Recall the definition of Hecke operators on modular symbols. Formally endow the set of maps from $\deltazero$ to $\D_{k,k}(\roi_L)$ with the action of an operator $U_p$ defined as
\[U_p \defeq\displaystyle\sum\limits_{[\alpha] \in \roi_K/(p)}\matrd{1}{\alpha}{0}{p}.\]

Here $\gamma = \smallmatrd{a}{b}{c}{d}$ acts as 
\[f|\gamma(D) = f(\gamma D)|\gamma,\]
and we consider such matrices as acting via the embedding of $\GLt(\roi_K) \hookrightarrow \GLt(\roi_L)\times\GLt(\roi_L)$, arising from equation (\ref{embedding}). Note that the image of $\smallmatrd{1}{\alpha}{0}{p}$ under this embedding has the form as described in the Lemma.
\begin{mrem}
Note that we've made some choice of orbit representatives for the action of $\Gamma$ on the double coset $\Gamma\smallmatrd{1}{0}{0}{p}\Gamma$. If we consider $U_p$ as a double coset operator on \emph{set-theoretic} maps, this is not well-defined up to choice of such representatives. However, as long as we are consistent in our choice, it shall not matter; hence we simply define $U_p$ in this very specific way and ignore where it comes from.
\end{mrem}

\subsection{Summary}\label{reformulation}
Before embarking on the proof of a suitable lifting theorem, we first take stock of the work so far. We have a monoid $\Sigma_0(\roi_L)$ acting on an $\roi_L$-module $\D_{k,k}(\roi_L)$, with  a $\Sigma_0(\roi_L)$-stable filtration $F^N\D_{k,k}(\roi_L)$ of $\D_{k,k}(\roi_L)$ leading to an inverse system $(A^N\D_{k,k}(\roi_L))$ of $\roi_L$-modules satisfying 
\[\lim\limits_{\longleftarrow}A^N\D_{k,k}(\roi_L) = \D_{k,k}(\roi_L).\]
Furthermore, $A^0\D_{k,k}(\roi_L) \cong V_{k,k}^*(\roi_L)$. Take $\lambda \in L^*$ with $v(\lambda) < k+1$, and define 
\[D = \{\mu \in \D_{k,k}(\roi_L): \rho^0(\mu) \in V_{k,k}^\lambda(\roi_L)\}.\]
Suppose $\gamma$ is a summand of the $U_p$ operator defined above; then Lemma \ref{tech}(i) tells us that if $\mu \in D$, then $\mu|\gamma \in \lambda \D_{k,k}(\roi_L)$, while Lemma \ref{tech}(ii) says that if $\mu \in F^N\D_{k,k}(\roi_L)$ then $\mu|\gamma \in \lambda F^{N+1}\D_{k,k}(\roi_L).$\\
\\
This then gives us the exact situation described in the next section.


\section{A Lifting Theorem}
The following is an abstraction of some of the elements of Greenberg's work in \cite{Gre07}. Since we will use essentially the same ideas multiple times, it is presented here in a completely general setting. The notation given is suggestive, and throughout, the reader should imagine the objects of (i) to (vii) below to be the obvious analogues from Section \ref{reformulation}.\\
\\
Suppose $L/\Qp$ is a finite extension, and that we have:
\begin{itemize}
\item[(i)] a monoid $\Sigma$,
\item[(ii)] a $\roi_L$-module $D$ that has a right action of $\Sigma$,
\item[(iii)] a $\Sigma$-stable filtration of $D$, $D \supset \mathcal{F}^0D \supset  \f^1D \supset \cdots,$ where if we define $\mathcal{A}^ND \defeq D/\mathcal{F}^ND$, then we have
\[\lim\limits_{\longleftarrow}\mathcal{A}^ND = D,\]
and where the $\mathcal{F}^ND$ have trivial intersection,
\item[(iv)] a right $\Sigma$-stable submodule $A$ of $\mathcal{A}^0D$, denoting $D_A \defeq \{\mu \in D: \mu \newmod{\mathcal{F}^0D} \in A\},$
\item[(v)] an operator $U = \sum\limits_{i=0}^r \gamma_i,$ where $\gamma_i \in \Sigma$,
\item[(vi)] a subgroup $\Gamma \leq \Sigma$ such that for each $j$, we have
\[\Gamma \gamma_j\Gamma = \coprod\limits_{i=0}^r\Gamma\gamma_i,\]
\item[(vii)] and a (countable) left $\Z[\Gamma]$-module $\Delta$.
\end{itemize}
For a right $\Z[\Sigma]$-module $\D$, endow the space of homomorphisms from $\Delta$ to $\D$ with a right $\Sigma$-action by
\[(\phi|\gamma)(E) = \phi(\gamma\cdot E)|\gamma.\]
For such $\D$, write $\symb_\Gamma(\D) = \Hom_\Gamma(\Delta,\D)$ for the space of $\Gamma$-equivariant homomophisms. Note that $U$ acts on this space.
\begin{mthm}\label{liftingtheorem}
Suppose that $\lambda$ is a non-zero element of $\roi_L$, that $D_A$ and $A$ have trivial $\lambda$-torsion, and that for each $\gamma_i$ appearing in the $U$ operator, we have:
\begin{itemize}
\item[(a)] If $\mu \in D_A$, then $\mu|\gamma_i \in \lambda D,$ and
\item[(b)]If $\mu \in \mathcal{F}^ND$, then $\mu|\gamma_i \in \lambda\mathcal{F}^{N+1}D$.
\end{itemize}
Then there is a natural map $\rho^0:\symb_\Gamma(D_A) \rightarrow \symb_\Gamma(A)$ whose restriction to the $\lambda$-eigenspaces of the $U$ operator is an isomorphism.
\end{mthm}
For clarity, the proof will be broken into a series of smaller steps. We have natural $\Sigma$-equivariant projection maps 
\[\pi^N: D \longrightarrow \mathcal{A}^ND\]that induce $\Sigma$-equivariant maps 
\[\rho^N: \symb_\Gamma(D) \longrightarrow \symb_\Gamma(\mathcal{A}^ND),\]
(and hence $\rho^0: \symb_\Gamma(D_A) \rightarrow \symb_\Gamma(A)$ by restriction) as well as maps $\pi^{M,N}:\mathcal{A}^MD \rightarrow \mathcal{A}^ND$ for $M \geq N$ that similarly induce maps $\rho^{M,N}.$ Thus we have an inverse system, and also it's straightforward to see that
\begin{align*}
\lim\limits_{\longleftarrow}\symb_\Gamma(\mathcal{A}^ND) = \symb_\Gamma(D).
\end{align*}

First we pass to a filtration where the $\Sigma$-action is nicer. Define $\mathcal{F}^ND_A = \mathcal{F}^ND \cap D_A.$ This is a $\Sigma$-stable filtration of $D_A$, since $A$ is $\Sigma$-stable and the projection maps are $\Sigma$-equivariant. Define  $\mathcal{A}^ND_A = D_A/\mathcal{F}^ND_A$, so that we have the following (where the vertical maps are injections):
\begin{diagram}
 D&&\rTo^{\pi^M}&&\mathcal{A}^MD&&\rTo^{\pi^{M,N}}&&\mathcal{A}^ND \\
\uTo&&&&\uTo&&&&\uTo \\
 D_A&&\rTo^{\pi^M}&&\mathcal{A}^MD_A&&\rTo^{\pi^{M,N}}&&\mathcal{A}^ND_A
\end{diagram}
Again, we see easily that
\begin{align}
\label{ila}\lim\limits_{\longleftarrow}\symb_\Gamma(\mathcal{A}^ND_A) = \symb_\Gamma(D_A).
\end{align}
Firstly, since $\mathcal{A}^ND_A$ may have non-trivial $\lambda$-torsion, we should make the statement ``$U$-eigensymbol in $\symb_\Gamma(\mathcal{A}^ND_A)$'' more precise. By condition (b) of \ref{liftingtheorem}, if $\gamma$ is a summand of $U$, and $\mu \in D_A$, then $\mu|\gamma \in \lambda D$. Accordingly, given a $\Z$-homomorphism $\varphi$ from $\Delta$ to $D_A$, we have $(\varphi|\gamma)(E) = \varphi(\gamma E)|\gamma = \lambda x$, for $E \in \Delta$ and some $x \in D$. Define a formal operator 
\[V_\gamma: \Hom(\Delta,D_A) \longrightarrow \Hom(\Delta, D) \]
by
\[(\varphi|V_\gamma)(E) = x,\]
so that we have an equality of operators $\lambda V_\gamma = \gamma|_{\Hom(\Delta,D_A)}$. 
\begin{mrems}\begin{itemize}
\item[(i)]Note that as $\rho_0$ is $\Sigma$-equivariant and $A$ is $\Sigma$-stable, $D_A$ is $\Sigma$-stable, so $\gamma|_{\Hom(\Delta,D_A)}$ is indeed an operator on $\Hom(\Delta,D_A)$. 
\item[(ii)]The reason we don't simply just define $V_\gamma = \lambda^{-1}\gamma$ is that `dividing by $\lambda$' is not in general a well-defined notion on $D$. 
\end{itemize}
\end{mrems}
Further define 
\[V = \sum\limits_{i=0}^r V_{\gamma_i},\]
so that we have an equality of operators $\lambda V = U|_{\Hom(\Delta,D_A)}$ (where the right hand side is an operator on $\Hom(\Delta,D_A)$ by Remark (i) above).\\
\\
It is easy to see that, for each $N$, $V$ gives rise to an operator $V_N : \mathcal{A}^ND_A \rightarrow \mathcal{A}^ND$. Under the natural injection $\mathcal{A}^ND_A \hookrightarrow \mathcal{A}^ND$, we may consider $\symb_\Gamma(\mathcal{A}^ND_A)$ as a subset of $\symb_\Gamma(\mathcal{A}^ND)$. We say an element $\varphi_N \in \symb_\Gamma(\mathcal{A}^ND_A)$ is a \emph{$U$-eigensymbol of eigenvalue $\lambda$} if $\varphi_N|V_N$ lives in $\symb_\Gamma(\mathcal{A}^ND_A)$ and moreover we have $\varphi_N|V_N = \varphi_N$.\\
\\
Take a $U$-eigensymbol $\phi_0 \in \symb_\Gamma(A) = \symb_\Gamma(\mathcal{A}^0D_A)$ with eigenvalue $\lambda$. Suppose a lift to a $U$-eigensymbol $\phi_N \in \symb_\Gamma(\mathcal{A}^ND_A)$ exists. We can take an arbitrary lift of $\phi_N$ to some $\Z$-homomorphism
\[\phi:\Delta \longrightarrow D.\]
Such a lift exists, as we can take some $\Z$-basis of $\Delta$ (using countability) and define $\phi$ on this basis, extending $\Z$-linearly (noting that this gives a well-defined lift since $\phi_N$ is also a $\Z$-homomorphism).\\
\\
Now, since $\phi_N$ is a $U$-eigensymbol with eigenvalue $\lambda$, it follows that $\phi|V$ is also a lift of $\phi_N$ to $\Hom(\Delta,D_A)$. The maps $\pi^N$, inducing the maps $\rho^N$, can be used immediately to extend the definition of $\rho^N$ to the space of $\Z$-homomorphisms from $\Delta$ to $D_A$ (rather than just the ones that are $\Gamma$-equivariant). Define
\[\phi_{N+1} = \rho^{N+1}\left(\phi|V\right),\]
a $\Z$-homomorphism from $\Delta$ to $\mathcal{A}^{N+1}D_A$. Note that since the maps $\rho^{N}, \rho^{M,N}$ are $\Sigma$-equivariant (and hence $V$-equivariant), we have compatibility relations 
\[\rho^{N+1,N}(\phi_{N+1}) = \phi_{N},\]
so that the following lemma says that the family $\{\phi_N\}$ we obtain gives an element of the inverse limit given in equation (\ref{ila}).
\begin{mlem}\label{esymbol}
$\phi_{N+1}$ is a well-defined $U$-eigensymbol in $\symb_\Gamma(\mathcal{A}^{N+1}D_A).$
\end{mlem}
We prove this lemma in a series of claims.
\begin{mcla}\label{claim}
If $\gamma = \gamma_i$ is a summand of the $U$ operator, and $\phi'$ is another lift of $\phi_N$, then 
\[\rho_{N+1}\left(\phi|V_\gamma\right) = \rho_{N+1}\left(\phi'|V_\gamma\right).\]
In particular, $\phi_{N+1}$ is independent of the choice of $\phi$ above $\phi_N$.
\end{mcla}
\begin{proof}
To say that $\phi$ and $\phi'$ are both lifts of $\phi_N$ is to say that the image of $\phi - \phi'$ under $\rho^N$ takes the value 0 in $\mathcal{A}^ND_A,$ that is, $\phi-\phi'$ takes values in $\mathcal{F}^ND_A.$ Thus by condition (b) in Theorem \ref{liftingtheorem}, 
\[(\phi-\phi')\bigg|\gamma \in \lambda \mathcal{F}^{N+1}D,\]
that is,
\[(\phi - \phi')\bigg|V_\gamma \in \mathcal{F}^{N+1}D.\]
As $\phi_N$ is a $U$-eigensymbol of eigenvalue $\lambda$, the image of $(\phi - \phi')|V_\gamma$ under $\rho^{N+1}$ lies in $\mathcal{A}^{N+1}D_A$, and by the above, it is 0. But this is precisely the statement that we required.
\end{proof}

\begin{mcla}$\phi_{N+1}$ is $\Gamma$-equivariant.
\end{mcla}
\begin{proof}
Let $\gamma \in \Gamma.$ As the map $\rho^{N+1}$ is $\Sigma$-equivariant, it follows that
\[\phi_{N+1}\bigg|\gamma = \frac{1}{\lambda}\rho^{N+1}\left(\sum\limits_{i=0}^r\phi\bigg|_k(\gamma_i\gamma)\right),\]
where the division by $\lambda$ is purely formal and well-defined by the remarks above. By condition (vi) above, we have a double coset decomposition 
\[\Gamma\gamma_{j}\Gamma = \coprod\limits_{i=0}^{r}\Gamma\gamma_{i},\]
hence we can find $\gamma^j\in \Gamma$ such that
\[\sum\limits_{i=0}^r\phi\bigg|_k(\gamma_{i}\gamma) = \sum\limits_{j=0}^r\left(\phi\bigg|_k\gamma^{j}\right)\bigg|_k\gamma_{j}.\]
Since $\phi_{N}$ is $\Gamma$-invariant, $\phi|\gamma^{j}$ is a lift of $\phi_N$, and hence by Claim \ref{claim} it follows that
\begin{align*}
\phi_{N+1}\bigg|\gamma &= \frac{1}{\lambda}\rho^{N+1}\left(\sum\limits_{j=0}^r\left(\phi\bigg|_k\gamma^{j}\right)\bigg|\gamma_{j}\right)\\
&= \frac{1}{\lambda}\rho^{N+1}\left( \sum\limits_{j=0}^r\phi\bigg|\gamma_{j}\right) = \phi_{N+1},
\end{align*}
as required.
\end{proof}
\begin{mcla}$\phi_{N+1}$ is a $U$-eigensymbol with eigenvalue $\lambda$.
\end{mcla}
\begin{proof}
Recall that when we consider the operator $V$ acting on $\mathcal{A}^ND_A$, we denote it $V_N$, and that $\phi_N$ is a $U$-eigensymbol with eigenvalue $\lambda$ if $\phi_N$ is a fixed point of $V_N$. Note then that $\phi|V$ is also a lift of $\phi_N$ to $\Hom(\Delta,D_A)$. In particular, $\phi|V$ also lives in $\Hom(\Delta,D_A)$, so we can apply $V$ to it again. Thus, by $\Sigma$-equivariance and Claim \ref{claim},
\begin{align*}\phi_{N+1}\bigg|V_{N+1} &= \rho^{N+1}\left(\phi\bigg|V\right)\bigg|V_{N+1}\\
&= \rho^{N+1}\left(\phi\bigg|V^2\right) = \phi^{N+1},
\end{align*}
as required.
\end{proof}
\begin{proof}(\emph{Theorem \ref{liftingtheorem}}).
Surjectivity follows from the results above; take an element \[\phi_0 \in \symb_\Gamma(A).\] Then for each $N$ we can construct $\phi_N \in \symb_\Gamma(\mathcal{A}^ND_A)$, compatibly with the projection maps $\pi^{M,N}$, and thus obtain a well-defined element of the inverse limit, which is the domain. By construction, this element has image $\phi_0$ under $\rho^0$.\\
\\
To prove injectivity, take some
\[\varphi \in \ker(\rho^0) = \symb_\Gamma(D_A\cap \mathcal{F}^0D).\]
Applying the operator $V$ recursively gives
\[\varphi = \varphi\bigg|V^N \in \mathcal{F}^ND_A;\]
indeed, by condition (b) in Theorem \ref{liftingtheorem}, $\varphi|V^N \in \f^ND$, and since $\varphi$ is a $U$-eigensymbol with eigenvalue $\lambda$, we can make sense of $\varphi|V^N = \varphi \in D_A$. Thus the image of $\varphi$ lies in the intersection of all the $\mathcal{F}^ND_A,$ that is, $\varphi = 0$ (by condition (iii)) as required.\\
\\Thus the map is a bijection, and thus an isomorphism, as required.
\end{proof}

An immediate consequence of this is the following corollary (the control theorem for Bianchi modular symbols):
\begin{mcor}\label{contthm1}
\begin{itemize}\item[(i)]Let $K/\Q$ be an imaginary quadratic field, $p$ a rational prime, and $L/\Qp$ the finite extension defined in Notation \ref{not}. Let $\Gamma = \Gamma_i$ be a twist of $\Gamma_1(\n) \leq \SLt(\roi_K)$ with $(p)|\n.$\\
\\
Let $\lambda \in L^\times$. Then, when $v_p(\lambda) < k+1$, the restriction of the specialisation map
\[\rho^0:\symb_\Gamma(\D_{k,k}(L))^{U_p=\lambda} \longrightarrow \symb_\Gamma(V_{k,k}^*(L))^{U_p=\lambda}\]
(where the superscript $(U_p=\lambda)$ denotes the $\lambda$-eigenspace for $U_p$) is an isomorphism.
\item[(ii)]In the same set up, with $\Omega_1(\n)$ defined as in Definition \ref{omega}, we have an isomorphism
\[\rho^0:\symb_{\Omega_1(\n)}(\D_{k,k}(L))^{U_p=\lambda} \longrightarrow \symb_{\Omega_1(\n)}(V_{k,k}^*(L))^{U_p=\lambda}.\]
\end{itemize}
\end{mcor}
\begin{proof}
To prove (i), recall that in the set up of Section \ref{reformulation}, Theorem \ref{liftingtheorem} says that the restriction of $\rho^0$ to the map 
\[\rho^0: \symb_\Gamma(D)^{U_p = \lambda} \longrightarrow \symb_\Gamma(V_{k,k}^\lambda(\roi_L))^{U_p = \lambda}\]
is an isomorphism. The result now follows by right-exactness of tensor product and Proposition \ref{tensoring}, since $D \otimes_{\roi_L} L \cong \D_{k,k}(L)$ and $V_{k,k}^\lambda \otimes_{\roi_L}K \cong V_{k,k}^*(L)$. Part (ii) is a trivial consequence as the $U_p$ operator acts separately on each component.
\end{proof}


\section{Values of Overconvergent Lifts}
This section will examine the spaces in which overconvergent lifts take values, refining our earlier results. Recall that $\Gamma = \Gamma_i$ is a twist of $\Gamma_1(\n)$.

\subsection{Locally Analytic Distributions}
We've shown that any classical Bianchi eigensymbol of suitable slope can be lifted to an overconvergent symbol that takes values in a space of rigid analytic distributions. However, the space we're truly interested in is a smaller space of distributions. A $p$-adic $L$-function should be a function on characters; but a rigid analytic distribution can take as input only functions that can be written as a convergent power series, and most finite order characters are locally constant, and thus cannot be written in this form. Instead, we want our lift to take values in the dual of \emph{locally analytic} functions.
\begin{mdef}\label{rsnbhd} Let $r,s \in \mathbb{R}_{>0}.$ Define the \emph{$(r,s)-$neighbourhood of $\roi_K\otimes_{\Z}\Zp$ in $\mathbb{C}_p^2$} to be
\[B(\roi_K\otimes_{\Z}\Zp,r,s) = \{(x,y)\in \Cp^2: \exists u \in \roi_K\otimes_{\Z}\Zp \text{ such that }|x-\sigma_1(u)|\leq r, |y-\sigma_2(u)|\leq s\},\]
where $\sigma_1,\sigma_2$ are the projections of the embedding $\sigma:K\otimes \Qp \hookrightarrow L^2$ to the first and second factors respectively.
\end{mdef}
When $p$ is split and $r$ and $s$ are both at least 1, $B(\roi_K\otimes_{\Z}\Zp,r,s)$ is the cartesian product of the closed discs of radii $r$ and $s$ in $\Cp$; and if $r = s = 1/p$, this is the cartesian product of two copies of the disjoint union of closed discs of radius $1/p$ with centres at 0, 1, ..., $p-1$.
\begin{mdef}Let $r$ and $s$ be as above. Then define the space of \emph{$L$-valued locally analytic functions of order $(r,s)$ on $\roi_K\otimes_{\Z}\Zp$}, denoted $\A[L,r,s]$, to be the space of rigid analytic functions on $B(\roi_K\otimes_{\Z}\Zp,r,s)$ that are defined over $L$.
\end{mdef}
\begin{mdef}Define the space of \emph{$L$-valued locally analytic distributions of order ($r,s$) on $\roi_K\otimes_{\Z}\Zp$} to be
\[\D[L,r,s] = \Hom_{\text{cts}}(\A[L,r,s],L).\]
\end{mdef}
We endow $\A[L,r,s]$ with a weight $k$ action of $\Sigma_0(\roi_K\otimes_{\Z}\Zp)$ identical to the action defined earlier on $\A_2(L)$; it is obvious that this action extends immediately to the larger space. It's then clear that by dualising, the action we obtain on $\D[L,r,s]$ is the restriction of the action on $\D_2(L)$. When talking about these spaces equipped with these actions, we denote them $\A_{k,k}[L,r,s]$ and $\D_{k,k}[L,r,s]$.\\
\\
For $r\leq r'$ and $s\leq s',$ we have a natural and completely continuous injection $\A[L,r',s']\hookrightarrow \A_{k,k}[L,r,s]$, since $B(\roi_K\otimes_{\Z}\Zp,r,s) \subset B(\roi_K\otimes_{\Z}\Zp,r',s')$. Since the polynomials are dense in each of these spaces, the image of this injection is dense. Using this compatibility, we make the following definitions:
\begin{mdef}Define the space of \emph{$L$-valued locally analytic functions on $\roi_K\otimes_{\Z}\Zp$} to be the direct limit
\[\mathcal{A}_{k,k}(L) \defeq \lim\limits_{\longrightarrow}\A_{k,k}[L,r,s] = \bigcup\limits_{r,s}\A[L,r,s].\]
\end{mdef}
\begin{mdef}Define the space of \emph{$L$-valued locally analytic distributions on $\roi_K\otimes_{\Z}\Zp$} to be
\[\mathcal{D}_{k,k}(L) \defeq \Hom_{\text{cts}}(\mathcal{A}_{k,k}(L),L).\]
\end{mdef}
\begin{mprop}\label{invlimiso}There is a canonical $\Sigma_0(\roi_K\otimes_{\Z}\Zp)$-equivariant isomorphism
\[\mathcal{D}_{k,k}(L) \cong \lim\limits_{\longleftarrow}\D_{k,k}[L,r,s] = \bigcap\limits_{r,s}\D_{k,k}[L,r,s].\]
\end{mprop}
\begin{proof}
The natural maps in each direction are inverse to each other and $\Sigma_0(\roi_K\otimes_{\Z}\Zp)$-equivariant.
\end{proof}

\subsection{The Action of $\Sigma_0(\roi_K\otimes_{\Z}\Zp)$}\label{actionofsigma}
The action of certain elements of $\Sigma_0(\roi_K\otimes_{\Z}\Zp)$ naturally moves us up and down the direct/inverse systems.

\begin{mlem}\label{overconverge}
\begin{itemize}
\item[(i)] Let $g \in \A_{k,k}[L,r,s]$ and $a \in \roi_K\otimes_{\Z}\Zp.$ Then
\[\gamma\cdot_{k}g,\hspace{12pt}\gamma = \matrd{1}{a}{0}{p^n},\]
naturally extends to $B(\roi_K\otimes_{\Z}\Zp,rp^n,sp^n)$ and thus gives an element of $\A_{k,k}[L,rp^{n},sp^{n}].$
\item[(ii)] Let $\mu \in \D_{k,k}[L,r,s],$ and $\gamma$ as above. Then $\mu|_{k}\gamma$ naturally gives an element of the smaller space $\D_{k,k}[L,rp^{-n},sp^{-n}].$
\end{itemize}
\end{mlem}
\begin{proof}
For $(x,y) \in B(\roi_K\otimes_{\Z}\Zp,rp^n,sp^n)$, there exists $b \in \roi_K\otimes_{\Z}\Zp$ such that $|x-\sigma_1(b)| \leq rp^n$ and $|y-\sigma_2(b)| \leq sp^n$, where $\sigma_1,\sigma_2$ are as defined after Definition \ref{rsnbhd}. Then $|(\sigma_1(a) + p^nx ) - (\sigma_1(a) + p^n\sigma_1(b))| \leq r,$ and similarly for $y$, so that 
\[(\sigma_1(a) + p^nx, \sigma_2(a) + p^ny) \in B(\roi_K\otimes_{\Z}\Zp,r,s).\]
For such $x,y$, we have
\[\gamma\cdot_k g(x,y) = g(\sigma_1(a) + p^nx,\sigma_2(a) + p^ny),\]
and since $g$ is defined on $B(\roi_K\otimes_{\Z}\Zp,r,s)$, the result follows.\\
\\
For part (ii), note that the action of $\gamma$ gives a map 
\[\A_{k,k}[L,r,s] \longrightarrow \A_{k,k}[L,rp^n,sp^n],\]
and hence dualising, the action gives a map
\[\D_{k,k}[L,rp^n,sp^n] \longrightarrow \D_{k,k}[L,r,s].\]
This is the required result (though scaled by a factor of $p^n$).
\end{proof}
\begin{mprop}\label{locanalyticprop}
Suppose that $\Psi \in \symb_\Gamma(\D_{k,k}(L))$ is a $U_p$-eigensymbol with non-zero eigenvalue. Then $\Psi$ is an element of $\symb_\Gamma(\mathcal{D}_{k,k}(L)).$
\end{mprop}
\begin{proof}
Firstly, $\D_{k,k}(L) = \D_{k,k}[L,1,1].$ Note that $U_p$ acts invertibly on the $U_p$-eigenspace, so that for each integer $n$, there exists some eigensymbol $\Psi'$ with $\Psi = \Psi'|U_p^n.$ The $U_p^n$ operator can be described explicitly as
\[U_p^n = \sum\limits_{[a]\in \roi_K/(p^n)}\matrd{1}{a}{0}{p^n},\]
so combining with Lemma \ref{overconverge} shows that $\Psi$ takes values in $\D_{k,k}[L,p^{-n},p^{-n}]$ for each $n$, and thus in $\lim\limits_{\longleftarrow}\D_{k,k}[L,r,s] \cong \mathcal{D}_{k,k}(L)$, using Proposition \ref{invlimiso}. The result follows.
\end{proof}
\begin{mcor}\begin{itemize}
\item[(i)]Let $K/\Q$ be an imaginary quadratic field, $p$ a rational prime, and $L/\Qp$ the finite extension defined in Notation \ref{not}. Let $\Gamma = \Gamma_i$ be a twist of $\Gamma_1(\mathfrak{n}) \leq \SLt(\roi_K)$ with $(p)|\mathfrak{n}.$\\
\\
Let $\lambda \in L^\times$. Then, when $v(\lambda) < k+1,$ the restriction of the natural map
\[\symb_\Gamma(\mathcal{D}_{k,k}(L))^{U_p = \lambda} \longrightarrow \symb_\Gamma(V_{k,k}^*(L))^{U_p=\lambda}\]
is an isomorphism.
\item[(ii)]In the set-up of above, and with $\Omega_1(\n)$ as defined in Definition \ref{omega}, the restriction of the natural map
\[\symb_{\Omega_1(\n)}(\mathcal{D}_{k,k}(L))^{U_p = \lambda} \longrightarrow \symb_{\Omega_1(\n)}(V_{k,k}^*(L))^{U_p=\lambda}\]
is an isomorphism.
\end{itemize}
\end{mcor}

\subsection{Admissible Distributions}
For each pair $r,s,$ the space $\D_{k,k}[L,r,s]$ admits a natural operator norm $||\cdot ||_{r,s}$ via
\[||\mu||_{r,s} = \sup\limits_{0\neq f\in \A_{k,k}[L,r,s]}\frac{|\mu(f)|_p}{|f|_{r,s}},\]
where $|\cdot |_p$ is the usual $p$-adic absolute value on $L$ and $|\cdot |_{r,s}$ is the sup norm on $\A_{k,k}[L,r,s].$ Note that if $r \leq r', s\leq s'$, then $||\mu||_{r,s} \geq ||\mu||_{r',s'}$ for $\mu \in \D_{k,k}[L,r',s'].$\\
\\
These norms give rise to a family of norms on the space of locally analytic functions. It is natural to classify locally analytic distributions by growth properties as we vary in this family. 
\begin{mdef}\label{admissible}Let  $\mu \in \mathcal{D}_{k,k}(L)$ be a locally analytic distribution.
We say $\mu$ is \emph{$h$-admissible} if 
\[ ||\mu||_{r,r} = O(r^{-h})\]
as $r \rightarrow 0^+.$
\end{mdef}
The following lemma is a useful technical result describing the family of norms of a $\Gamma$-orbit in $\A_{k,k}[L,r,s],$ where $\Gamma$ is as before. It gives universal constants that will be useful in the sequel. 
\begin{mlem}\label{equivnorm}
There exist positive constants $C$ and $C'$ such that
\[C|\gamma\cdot_{k}f|_{r,s} \leq |f|_{r,s} \leq C'|\gamma\cdot_{k}f|_{r,s}\]
for every $\gamma \in \Gamma$ and $f \in \A_{k,k}[L,r,s]$.
\end{mlem}
\begin{proof}
The action of $\gamma$ by
\[\gamma\cdot(x,y) = \left(\frac{b + dx}{a+cx},\frac{b'+d'y}{a' + c'y}\right),\hsp \sigma_1(\gamma) = \matr, \sigma_2(\gamma) = \matrd{a'}{b'}{c'}{d'}\]
 maps $B(\roi_K\otimes_{\Z}\Zp,r,s)$ bijectively to itself. Furthermore, we have 
\[|(a+cx)^k(a'+c'y)^k|_p \leq \max\{1,|x|_p^k, |y|_p^k, |xy|_p^k\} \leq \max\{1,r^k, s^k, (rs)^k\} = C^{-1},\]
say, noting that $C^{-1}$ is certainly positive. Thus
\begin{align*}
|\gamma\cdot_{k}f|_{r,s} = &\sup\limits_{(x,y) \in B(\roi_K\otimes_{\Z}\Zp,r,s)}\bigg|(a+cx)^k(a'+c'y)^k\bigg|_p\cdot\bigg|f\bigg(\gamma\cdot(x,y)\bigg)\bigg|_p\\
\leq C^{-1}&\sup\limits_{(x,y) \in B(\roi_K\otimes_{\Z}\Zp,r,s)}|f(x,y)|_p\\
=C^{-1}&|f|_{r,s},
\end{align*}
from which the first inequality follows. The reverse direction follows from symmetry by considering the action of $\gamma^{-1}$.
\end{proof}
\begin{mdefprop}\label{defnorm} Let $\Psi \in \symb_\Gamma(\mathcal{D}_{k,k}(L)),$ and $r,s \leq 1.$ Defining 
\[||\Psi||_{r,s} \defeq \sup\limits_{D \in \Delta_0}||\Phi(D)||_{r,s}\]
gives a well-defined norm on $\symb_\Gamma(\mathcal{D}_{k,k}(L))$.
\end{mdefprop}
\begin{proof}
Pick a finite set of generators $D_1, ... , D_m$ for $\Delta_0$ as a $\Z[\Gamma]-$module. Writing $D = \alpha_1 D_1 + \cdots \alpha_n D_n$ with $\alpha_i \in \Z[\Gamma]$, and using Lemma \ref{equivnorm}, we see that there exists a $C$ such that (without loss of generality)
\[||\Psi(D)||_{r,s} \leq C ||\Psi(D_1)||_{r,s}.\]
In particular, the supremum is finite and hence gives a well-defined norm, as required.
\end{proof}
The values of an overconvergent eigensymbol satisfy further conditions of the type above depending on their slope.
\begin{mprop}\label{admissibility}
Suppose $\Psi \in \symb_\Gamma(\mathcal{D}_{k,k}(L))$ is a $U_p$-eigensymbol with eigenvalue $\lambda$ and slope $h = v(\lambda).$ Then, for every $D \in \Delta_0$, the distribution $\Psi(D)$ is $h$-admissible.
\end{mprop}
\begin{proof}
For any $r$ and a positive integer $n$, we have
\begin{align*}
||\Psi(D)||_{\frac{r}{p^{n}},\frac{r}{p^{n}}} &= |\lambda|^{-n}\bigg|\bigg|\left(\Psi\bigg|_{k}U_p^n\right)(D)\bigg|\bigg|_{\frac{r}{p^{n}},\frac{r}{p^{n}}}\\
&\leq|\lambda|^{-n}\max\limits_{[a]\in \roi_K/(p^n)}\bigg|\bigg|\Psi\left(\matrd{1}{a}{0}{p^n}D\right)\bigg|_{k}\sigma\left[\matrd{1}{a}{0}{p^n}\right]\bigg|\bigg|_{\frac{r}{p^{n}},\frac{r}{p^{n}}},
\end{align*}
where $\sigma$ is the embedding $\roi_K \hookrightarrow \roi_L\times\roi_L$ of equation (\ref{embedding}),
\begin{align*}&\leq|\lambda|^{-n}\max\limits_{[a]\in \roi_K/(p^n)}\bigg|\bigg|\Psi\left(\matrd{1}{a}{0}{p^n}D\right)\bigg|\bigg|_{r,r}\\
&\leq |\lambda|^{-n}||\Psi||_{r,r},
\end{align*}
the norm defined in Definition-Proposition \ref{defnorm}. Here the second to last inequality follows since for $\gamma = \smallmatrd{1}{a}{0}{p^n} \in \Sigma_0(\roi_K\otimes_{\Z}\Zp),$ we have, for any $\mu \in \D_{k,k}[L,r,s],$
\[\bigg|\bigg|\mu\bigg|_{k}\gamma\bigg|\bigg|_{\frac{r}{p^{n}},\frac{s}{p^{n}}} \leq ||\mu||_{r,s}.\]
This is simply because, for $f \in \A_{k,k}[L,rp^{-n},sp^{-n}],$
\[||f||_{rp^{-n},sp^{-n}} \geq ||\gamma\cdot_{k}f||_{r,s},\]
and hence 
\[\bigg|\bigg|\mu\bigg|_{k}\gamma\bigg|\bigg|_{\frac{r}{p^{n}},\frac{s}{p^{n}}} = \sup\limits_{f \in \A_{k,k}[L,\frac{r}{p^{n}},\frac{s}{p^{n}}]}\frac{|\mu(\gamma\cdot_{k}f)|_p}{||f||_{\frac{r}{p^{n}},\frac{s}{p^{n}}}}\]
\[\leq \sup\limits_{f \in \A_{k,k}[L,\frac{r}{p^{n}},\frac{s}{p^{n}}]}\frac{|\mu(\gamma\cdot_{k}f)|_p}{ ||\gamma\cdot_{k}f||_{r,s}}\]
\[ \leq  \sup\limits_{g \in \A_{k,k}[L,r,s]}\frac{|\mu(g)|_p}{ ||g||_{r,s}} = ||\mu||_{r,s}.\qedhere\]
\end{proof}


\section{The Case $p$ Split}
The results above hold for the $U_p$ operator regardless of the splitting behaviour of the prime $p$ in $K$. When $p$ is inert, this is the whole story, and if $p$ is ramified, then there are only slight modifications to make to obtain a finer result for the $U_{\pri}$ operator, where $(p) =\pri^2$ in $K$. In the case that $p$ splits in $K$ as $\pri\pribar$, however, we can obtain more subtle results. The crux of this section is that we can prove a control theorem for a `small slope' condition that encompasses far more possible eigensymbols. For example, suppose $E/K$ is a modular elliptic curve with supersingular reduction at both primes above $p$. The associated weight $(k,k) = (0,0)$ eigenform gives rise to a Bianchi eigensymbol $\phi \in \symb_\Gamma(V_{k,k}^*(L))$ with slope $1/2$ at $\pri$ and $1/2$ at $\pribar$. This symbol will have slope $1 = k+1$ under the $U_p$ operator, and hence we cannot lift it using the control theorem proved above. The results below will allow us to lift even symbols such as this.\\
\\
Many of the results and proofs closely mirror those of previous sections. We first prove a more refined control theorem, again by writing down a suitable filtration and using Theorem \ref{liftingtheorem}, then go on to prove an admissibility result for such lifts.\\
\\
Throughout, we assume that $p$ splits as $\pri\pribar$. Whilst ultimately we want to prove a control theorem for the full space of modular symbols for $\Omega_1(\n)$, it is simpler to instead work with a principal power of $\pri$ and look at each individual component of the direct sum separately, combining together at the end. 

\subsection{Lifting Simultaneous Eigensymbols of $U_\pri$ and $U_{\overline{\pri}}$}
\label{simesymbols}
The following results will show that it's possible to lift a classical Bianchi symbol to a space of Bianchi modular symbols that are overconvergent in one variable, and then again from this space to the space of fully overconvergent Bianchi modular symbols we considered previously.\\
\\ 
To do so, consider the space 
\[\DV(R) = \Hom_{\text{cts}}(\A_k(R)\otimes_R V_k(R),R),\]
with the appropriate action of $\Sigma_0(\roi_L)^2$ (where this makes sense) induced from the action on $\A_{k,k}$. This gives us
\[V_{k,k}(R) \subset \DV(R) \subset \D_{k,k}(R).\]
Now put $R = \roi_L,$ and recall the filtration in the one variable case from Definition \ref{onevarfilt}. We now define new filtrations to reflect lifting by one variable at a time.
\begin{mdef}
Define
\begin{itemize}
\item[(i)]$\begin{aligned}[t]\mathcal{F}^N\DV(\roi_L) &= \mathcal{F}^N\D_k(\roi_L)\otimes_{\roi_L}V_k^*(\roi_L)\\
& = \{\mu \in \DV(\roi_L): \mu(x^iy^j) \in \pi_L^{N-i}\roi_L \text{ for all }j\}.
\end{aligned}$

\item[(ii)]$\begin{aligned}[t]\DV^0(\roi_L) &= \ker(\DV(\roi_L) \rightarrow V_{k,k}^*(\roi_L))\\
&= \{\mu \in \DV(\roi_L): \mu(x^iy^j) = 0 \text{ for all } 0 \leq i \leq k\}.
\end{aligned}$
\item[(iii)]$F^N\DV(\roi_L) = \mathcal{F}^N\DV(\roi_L) \cap \DV^0(\roi_L).$
\end{itemize}
\end{mdef}

\begin{mdef}Define
\begin{itemize}
\item[(i)]$\begin{aligned}[t]\mathcal{F}_\pri^N\D_{k,k}(\roi_L) &= \D_k(\roi_L) \ctp_{\roi_L}\mathcal{F}^N(\roi_L)\\
&= \{\mu \in \D_{k,k}(\roi_L): \mu(x^iy^j) \in \pi_L^{N-j}\roi_L \text{ for all }i\}.
\end{aligned}$
\item[(ii)]$\begin{aligned}[t]\D_{k,k,\pri}^0(\roi_L) &= \ker(\D_{k,k}(\roi_L) \rightarrow \DV(\roi_L))\\
&= \{\mu \in \D_{k,k}(\roi_L): \mu(x^iy^j) = 0 \text{ for all } j\geq 0\}.
\end{aligned}$
\item[(iii)]$F^N_\pri\D_{k,k}(\roi_L) = \mathcal{F}^N_\pri\D_{k,k}(\roi_L) \cap \D_{k,k,\pri}^0(\roi_L).$
\end{itemize}
\end{mdef}
Further define
\[A^N\DV(\roi_L) = \DV(\roi_L)/F^N\DV(\roi_L),\]
\[A^N_\pri\D_{k,k}(\roi_L) = \D_{k,k}(\roi_L)/F^N_\pri\D_{k,k}(\roi_L).\]
Hence we now have filtrations
\begin{align*}
\mathcal{F}^0\DV(\roi_L) \subset \cdots\subset &\mathcal{F}^M\DV(\roi_L) \subset \cdots \subset \DV(\roi_L)\subset\\
 \cdots&\subset \mathcal{F}^N_\pri\D_{k,k}(\roi_L) \subset \cdots \subset \D_{k,k}(\roi_L).
\end{align*}
\begin{mprop}\label{twovarsigstab}
These filtrations are $\Sigma_0(\roi_L)^2$-stable.
\end{mprop}
\begin{proof}
This follows from the one variable case, as these filtrations are defined to be a tensor product of $\Sigma_0(\roi_L)$-stable spaces, and in the two variable case, $\Sigma_0(\roi_L)^2$ acts as $\Sigma_0(\roi_L)$ separately on each component.
\end{proof} 
These filtrations lead to $\Sigma_0(\roi_L)^2$-equivariant projection maps
\[\pi_1^N: \DV(\roi_L) \longrightarrow A^N\DV(\roi_L)\]
and
\[\pi_2^N: \D_{k,k}(\roi_L) \longrightarrow A_\pri^N\D_{k,k}(\roi_L),\]
which again give maps $\rho^N_1$ and $\rho^N_2$ on the corresponding symbol spaces.\\
\\
We want an analogue of Lemma \ref{tech} for this setting. Choose $n$ such that $\pri^n = (\beta)$ is principal (noting that this also forces $\overline{\pri}^n$ to be principal). Then we define $U_{\pri^n} = U_{\pri}^n$ as 
\[\sum\limits_{a\newmod{\pri^n}}\matrd{1}{a}{0}{\beta}.\]
We will prove control theorems for eigenspaces of the operators $U_{\pri}^n$ for $\pri|p$, which will give us the theorem for the operators $U_{\pri}$, as required. We can write down the appropriate analogue of $V_{k,k}^\lambda(\roi_L)$ in this setting as follows:

\begin{mdef}
Let $\lambda \in L^*$. Define
\begin{itemize}
\item[(i)]$\begin{aligned}[t]V_{k,k,\pri}^\lambda(\roi_L) &= V_k^\lambda(\roi_L) \otimes_{\roi_L} V_k^*(\roi_L)\\
&= \{f \in V_{k,k}^*(\roi_L): f(x^iy^j) \in \lambda p^{-i}\roi_L \text{ for }0 \leq i \leq \lfloor v(\lambda)\rfloor\},
\end{aligned}$
\item[(ii)]$\begin{aligned}[t]\DV&^\lambda(\roi_L) = \D_k(\roi_L) \otimes_{\roi_L}V_k^\lambda(\roi_L)\\ 
&= \{f \in \DV(\roi_L): f(x^iy^j) \in \lambda p^{-j}\roi_L \text{ for }0 \leq j \leq \lfloor v(\lambda) \rfloor\}.
\end{aligned}$
\end{itemize}
\end{mdef}
This gives the following situation:
\[\D_{k,k}(\roi_L) \labelrightarrow{\pi_2^0} \DV(\roi_L) \supset \DV^\lambda(\roi_L),\]
\[\DV(\roi_L) \labelrightarrow{\pi_1^0} V_{k,k}^*(\roi_L) \supset V_{k,k,\pri}^\lambda(\roi_L).\]
\begin{mprop}
These modules are $\Sigma_0(\roi_L)^2$-stable.
\end{mprop}
\begin{proof}
As in Proposition \ref{twovarsigstab}, this follows from the one variable case, as these are nothing but a tensor product of $\Sigma_0(\roi_L)$-stable spaces.
\end{proof}

The following two lemmas are practically identical in spirit to Lemma \ref{tech}, but they are included here for completeness. Recall that $n$ is a positive integer such that $\pri^n = (\beta)$ is principal.
\begin{mlem}\label{liftfirstvar}Let $a_1, a_2 \in \roi_L$.
\begin{itemize}
\item[(i)]Suppose $\mu \in \DV(\roi_L)$ with $\pi^0_1(\mu) \in V_{k,k,\pri}^\lambda(\roi_L)$. Then
\[\mu\bigg|_k\left[\matrd{1}{a_1}{0}{\beta}, \matrd{1}{a_2}{0}{\overline{\beta}}\right] \in \lambda\DV(\roi_L).\]
\item[(ii)] Suppose $v(\lambda) < n(k+1)$. Then for $\mu \in F^N\DV(\roi_L)$, we have 
\[\mu\bigg|_k\left[\matrd{1}{a_1}{0}{\beta}, \matrd{1}{a_2}{0}{\overline{\beta}}\right] \in \lambda F^{N+1}\DV(\roi_L).\]
\end{itemize}
\end{mlem}
\begin{proof}Take some $\mu \in \DV(\roi_L).$ Then
\[\mu\bigg|_k\left[\matrd{1}{a_1}{0}{\beta},\matrd{1}{a_2}{0}{\overline{\beta}}\right](x^my^n) = \mu((a_1+\beta x)^m(a_2+\overline{\beta}y)^n)\]
\[=\sum\limits_{i=0}^m\sum\limits_{j=0}^n\binomc{m}{i}\binomc{n}{j}a_1^{m-i}a_2^{n-j}\beta^{i}\overline{\beta}^j\mu(x^iy^j).\]
\begin{itemize}
\item[(i)] Suppose that $\pi^0(\mu)$ lies in $V_{k,k,\pri}^\lambda(\roi_L),$ so that $\mu(x^iy^j) \in \lambda p^{-i}\roi_L$ for any $i \leq\lfloor v(\lambda)\rfloor$. As $\beta^i \in p^{ni}\roi_L$, and $n\geq 1$, it follows that each term of the sum above lies in $\lambda\roi_L$, and hence we have the result. If instead $i$ is greater than $\left\lfloor v(\lambda)\right\rfloor$, it follows that $i > v(\lambda),$ so that $\beta^i \in \lambda\roi_L$, and hence the result follows as $\mu(x^iy^j) \in \roi_L$.
\item[(ii)] Now suppose $\mu \in F^N\DV(\roi_L)$. Again considering the sum above, the terms where $i \leq k$ vanish. If $i> k$, then 
\[ni \geq n(k+1) > v(\lambda),\]
since $\lambda$ has $p$-adic valuation $<n(k+1)$. As $\beta^{i}$ and $\lambda$ are divisible by integral powers of $\pi_L$, it follows that $\beta^{i} \in \pi_L\lambda\roi_L.$ Hence, as $\mu(x^iy^j) \in \pi_L^{N-i}\roi_L$, it follows that
\[\beta^i\mu(x^iy^j) \in \lambda \pi_L^{(N+1)-i}\roi_L,\]
which completes the proof.\qedhere
\end{itemize}
\end{proof}
Now let $\pribar^n = (\delta)$, with image $(\overline{\delta},\delta)$ in $\roi_L^2$ under $\sigma$. Note that $v(\delta) = n,$ whilst $\overline{\delta}$ is a unit in $\roi_L$.
\begin{mlem}\label{liftsecondvar}Let $a_1, a_2 \in \roi_L$.
\begin{itemize}
\item[(i)]Suppose $\mu \in \D_{k,k}(\roi_L)$ with $\pi^0_2(\mu) \in \DV^\lambda(\roi_L)$. Then
\[\mu\bigg|_k\left[\matrd{1}{a_1}{0}{\overline{\delta}}, \matrd{1}{a_2}{0}{\delta}\right] \in \lambda\D_{k,k}(\roi_L).\]
\item[(ii)] Suppose $v(\lambda) < n(k+1)$. Then for $\mu \in F^N_\pri\D_{k,k}(\roi_L)$, we have 
\[\mu\bigg|_k\left[\matrd{1}{a_1}{0}{\overline{\delta}}, \matrd{1}{a_2}{0}{\delta}\right] \in \lambda F^{N+1}_\pri\D_{k,k}(\roi_L).\]
\end{itemize}
\end{mlem}
\begin{proof}Identical to that of Lemma \ref{liftfirstvar} -- up to notation -- but with $j$'s replacing $i$'s where appropriate.
\end{proof}
We are hence in exactly the situation of Theorem \ref{liftingtheorem}, and applying it twice gives us:
\begin{mlem}\label{liftinghalveslemma}
Let $K/\Q$ be an imaginary quadratic field, $p$ a rational prime that splits as $\pri\overline{\pri}$ in $K$, $n$ an integer such that $\pri^n$ is principal, and $L/\Qp$ the finite extension defined in Notation \ref{not}. Let $\Gamma = \Gamma_i$ be a twist of $\Gamma_1(\n) \leq \SLt(\roi_K)$ with $(p)|\n.$ Let $\lambda \in L^\times$. Then, when $v(\lambda) < n(k+1)$, we have:
\begin{itemize}
\item[(i)] The restriction of the specialisation map
\[\rho_1^0:\symb_\Gamma(\DV(L))^{U_\pri^n=\lambda} \longrightarrow \symb_\Gamma(V_{k,k}^*(L))^{U_\pri^n=\lambda}\]
(where the superscript $(U_\pri=\lambda)$ denotes the $\lambda$-eigenspace for $U_\pri$) is an isomorphism.
\item[(ii)] The restriction of the specialisation map
\[\rho_2^0:\symb_\Gamma(\D_{k,k}(L))^{U_{\pribar}^n=\lambda} \longrightarrow \symb_\Gamma(\DV(L))^{U_{\pribar}^n=\lambda}\]
is an isomorphism.
\end{itemize}
\end{mlem}
\begin{proof}The only remaining details to fill in are formalities regarding tensor products, which are analogous to before and are omitted.
\end{proof}

We require results about the $U_\pri$ operator, whilst the results above are for the $U_{\pri}^n$ operator. We work around this using:
\begin{mlem}\label{untou}Suppose we have two $L$-vector spaces $D, V$ with a right action of an operator $U$, and a $U$-equivariant surjection $\rho: D \rightarrow V$ such that, for some positive integer $n$ and $\lambda \in L^\times$, the restriction of $\rho$ to the $\lambda^n$-eigenspaces of the $U^n$ operator is an isomorphism. Then the restriction of $\rho$ to the $\lambda$-eigenspaces of the $U$ operator is an isomorphism.
\end{mlem}
\begin{proof}
Suppose $\phi \in V$ is a $U$-eigensymbol with eigenvalue $\lambda$. Then $\phi$ is also a $U^n$-eigensymbol with eigenvalue $\lambda^n$, and accordingly there is a unique lift $\Psi$ of $\phi$ to a $U^n$-eigensymbol with eigenvalue $\lambda^n$. We claim that this is in fact a $U$-eigensymbol with eigenvalue $\lambda$. Indeed, by $U$-equivariance we have
\[\rho(\Psi|U) = \rho(\Psi)|U = \phi|U = \lambda\phi,\]
so that $\Psi|U$ is a lift of $\lambda\phi$. But clearly $\lambda\Psi$ is also a lift of $\lambda\phi$. Now as $\lambda\Psi$ and $\lambda\phi$ are both $U^n$-eigensymbols with eigenvalue $\lambda^n$, it follows that $\lambda\Psi$ is the unique eigenlift of $\lambda\phi$ under $\rho$; but then it follows that $\Psi|U = \lambda\Psi$, as required. The Lemma follows easily.
\end{proof}
This then gives:
\begin{mthm}\label{refinedcontrolthm}\begin{itemize}
\item[(i)]Let $K,$ $p,$ $\pri$, $n$, $\Gamma$ and $L$ be as in Lemma \ref{liftinghalveslemma}. Take $\lambda_1, \lambda_2 \in L^*$ with $v(\lambda_1), v(\lambda_2) < k+1.$ Then the restriction of the specialisation map
\[\rho^0:\symb_\Gamma(\D_{k,k}(L))^{U_\pri^n=\lambda_1^n, U_{\overline{\pri}}^n = \lambda_2^n} \longrightarrow \symb_\Gamma(V_{k,k}^*(L))^{U_\pri^n=\lambda_1^n, U_{\overline{\pri}}^n = \lambda_2^n}\]
(where the superscript denotes the simultaneous $\lambda_1^n$-eigenspace of $U_\pri^n$ and $\lambda_2^n$-eigenspace of $U_{\overline{\pri}}^n$) is an isomorphism.
\item[(ii)] In the set up of part (i), the restriction of the specialisation map 
\[\rho^0:\symb_{\Omega_1(\n)}(\D_{k,k}(L))^{U_\pri=\lambda_1, U_{\overline{\pri}} = \lambda_2} \longrightarrow \symb_{\Omega_1(\n)}(V_{k,k}^*(L))^{U_\pri=\lambda_1, U_{\overline{\pri}} = \lambda_2}\]
is an isomorphism.
\end{itemize}
\end{mthm}
\begin{proof}\begin{itemize}\item[(i)]
Take a simultaneous $\Upri^n$- and $\Upribar^n$-eigensymbol $\phi^0$, with eigenvalues $\lambda_1^n, \lambda_2^n$ respectively. Then Lemma \ref{liftinghalveslemma}(i) says that we can lift $\phi^0$ uniquely to some 
\[\varphi^0 \in \symb_\Gamma(\DV(L))^{\Upri^n = \lambda_1^n}.\]
We claim that $\varphi^0$ is a $\Upribar^n$-eigensymbol with eigenvalue $\lambda_2^n$. Indeed, consider the action of the operator $\lambda_2^{-n}\Upribar^n$. When applied to $\varphi^0$, the result is a $\Upri^n$-eigensymbol with eigenvalue $\lambda_1^n$, since
\begin{align*}\left(\varphi^0\bigg|\lambda_2^{-n}\Upribar^n\right)\bigg|\Upri^n &= \left(\varphi^0\bigg|\Upri^n\right)\bigg|\lambda_2^{-n}\Upribar^n\\
&= \lambda_1^n\varphi^0\bigg|\lambda_2^{-n}\Upribar^n,
\end{align*}
as $\Upri^n$ and $\Upribar^n$ commute. Then by the $\Sigma_0(\roi_L)^2$-equivariance of $\rho_1^0$ we have
\begin{align*}\rho_1^0\left(\varphi^0\bigg|\lambda_2^{-n}\Upribar^n\right) &= \rho_1^0(\varphi^0)\bigg|\lambda_2^{-n}\Upribar^n\\
&= \phi^0\bigg|\lambda_2^{-n}\Upribar^n = \phi^0,
\end{align*}
since $\phi^0$ is a $\Upribar^n$-eigensymbol with eigenvalue $\lambda_2^n$. But then by uniqueness, we must have
\[\varphi^0\bigg|\lambda_2^{-n}\Upribar^n = \varphi^0,\]
that is, $\varphi^0$ is a $\Upribar^n$-eigensymbol with eigenvalue $\lambda_2^n$, as required. Now we can use Lemma \ref{liftinghalveslemma}(ii) to lift $\varphi^0$ to some
\[\phi \in \symb_\Gamma(\D_{k,k}(L))^{\Upribar^n = \lambda_2^n}.\]
By an identical argument to that above, $\phi$ is a $\Upri^n$-eigensymbol with eigenvalue $\lambda_1^n$, and since by construction $\rho^0(\phi) = \phi^0$, this is the result.
\item[(ii)] From part (i), it is easy to see that we have an isomorphism
\[\rho^0:\symb_{\Omega_1(\n)}(\D_{k,k}(L))^{U_\pri^n=\lambda_1^n, U_{\overline{\pri}}^n = \lambda_2^n} \longrightarrow \symb_{\Omega_1(\n)}(V_{k,k}^*(L))^{U_\pri^n=\lambda_1^n, U_{\overline{\pri}}^n = \lambda_2^n}.\]
The result then follows directly from Lemma \ref{untou}, since there are well-defined $\Upri$ and $\Upribar$ operators on each of the spaces, and $\rho^0$ is equivariant with respect to these operators.
\end{itemize}
\end{proof}

\subsection{The Action of $\Sigma_0(\roi_K\otimes_{\Z}\Zp)$}
The following results are proved in an almost identical manner to those in Section \ref{actionofsigma}.

\begin{mlem}Let $p$ split as $\pri\overline{\pri}$ in $K$, with $\pri^n = (\beta)$ principal.
\begin{itemize}
\item[(i)] Let $g \in \A_{k,k}[L,r,s],$ and $a\in\roi_K\otimes_{\Z}\Zp.$ Then
\[\gamma\cdot_{k}g,\hspace{12pt}\gamma = \matrd{1}{a}{0}{\beta^m},\]
naturally extends to $B(\roi_K\otimes_{\Z}\Zp,rp^{mn},s)$ and thus gives an element of $\A_{k,k}[L,rp^{mn},s].$
\item[(ii)] Let $\mu \in \D_{k,k}[L,r,s],$ and $\gamma_i$ as above for $i = 1,2$. Then $\mu|_{k}\gamma$ naturally gives an element of the smaller space $\D_{k,k}[L,rp^{-mn},s].$
\end{itemize}
\end{mlem}
We also have an entirely analogous result for the $U_{\pribar}^n$ operator. Combining the two then gives the following:
\begin{mprop}
Suppose that $\Psi \in \symb_\Gamma(\D_{k,k}(L))$ is simultaneously a $U_{\pri}^n$- and $U_{\pribar}^n$-eigensymbol with non-zero eigenvalues. Then $\Psi$ is an element of $\symb_\Gamma(\mathcal{D}_{k,k}(L)).$
\end{mprop}

\begin{mcor}Let $K/\Q$ be an imaginary quadratic field, $p$ a rational prime that splits as $\pri\pribar$ in $K$, and $L/\Qp$ the finite extension defined in Notation \ref{not}. Let $\Omega_1(\n)$ be as defined in Definition \ref{omega}.\\
\\
Take $\lambda_1, \lambda_2 \in L^*$ with $v(\lambda_1), v(\lambda_2) < k+1.$ Then the restriction of the specialisation map
\[\rho^0:\symb_{\Omega_1(\n)}(\mathcal{D}_{k,k}(L))^{U_\pri=\lambda_1, U_{\overline{\pri}} = \lambda_2} \longrightarrow \symb_{\Omega_1(\n)}(V_{k,k}^*(L))^{U_\pri=\lambda_1, U_{\overline{\pri}} = \lambda_2}\]
is an isomorphism.
\end{mcor}

\subsection{Admissible Distributions}
In this new setting, we need a new definition of admissibility - namely one that encodes the slope at both $\pri$ and $\pribar$.
\begin{mdef}\label{splitadmissible}Let  $\mu \in \mathcal{D}_{k,k}(L)$ be a locally analytic distribution. We say $\mu$ is \emph{$(h_1,h_2)$-admissible} if
\[ ||\mu||_{r,s} = O(r^{-h_1})\]
uniformly in $s$ as $r\rightarrow 0^+$, and
\[ ||\mu||_{r,s} = O(s^{-h_2})\]
uniformly in $r$ as $s\rightarrow 0^+$.
\end{mdef}

\begin{mprop}
Suppose $p$ splits in $K$ as $\pri\overline{\pri}$, with $\pri^n = (\beta)$ principal, and suppose that $\Psi \in \symb_\Gamma(\mathcal{D}_{k,k}(L))$ is a $U_{\pri}^n$-eigensymbol with eigenvalue $\lambda_1^n$ and a $U_{\pribar}^n$-eigensymbol with eigenvalue $\lambda_2^n$ with slopes $h_i = v(\lambda_i).$ Then, for every $D \in \Delta_0$, the distribution $\Psi(D)$ is $(h_1,h_2)$-admissible.
\end{mprop}
\begin{proof}
The proof is similar to that of Proposition \ref{admissibility}; using $m$th powers of the operators $U_{\pri}^n$ and $U_{\pribar}^n$, we show that
\[||\Psi(D)||_{\frac{r}{p^{mn}},s} \leq |\lambda_1|^{-mn}||\Phi||_{r,s},\hsp \text{and} \hsp ||\Phi(D)||_{r,\frac{s}{p^{mn}}} \leq |\lambda_2|^{-mn}||\Phi||_{r,s}\]
respectively and combine to get the result.
\end{proof}

\section{The $p$-adic $L$-function}
In this final section, we combine the previous results to construct the $p$-adic $L$-function of a Bianchi modular form. Ultimately, we'll show that lifting the classical modular symbol attached to a small slope cuspidal Bianchi eigenform $\Phi$ under the control theorem and evaluating at $\{0\}-\{\infty\}$ gives us all of the information required to interpolate the critical values of its $L$-function, and gives us a sensible notion of this $p$-adic $L$-function.

\subsection{Evaluating at $\{0\}-\{\infty\}$}
\label{evalatzeroinf}
In the rational case, for a classical modular symbol $\phi_f$ attached to certain cusp forms $f$, and with overconvergent lift $\Psi_f$, the restriction of the distribution $\Psi_f(\{0\}-\{\infty\})$ to $\Zp^\times$ is the $p$-adic $L$-function of $f$. We want to emulate this result in the Bianchi case.
\begin{mnot}
To ease notation, we write $\roi_{K,p} \defeq \roi_K\otimes_{\Z}\Zp$.
\end{mnot}
\begin{mdef}Let $(\phi_1,...,\phi_h)$ be a classical Bianchi eigensymbol (resp. $\Phi$ a classical Bianchi eigenform), with $U_{\pri}$ eigenvalue(s) $a_{\pri}$ for $\pri|p$. We can canonically see $a_{\pri}$ as living in $\overline{\Qp}$ under the fixed embedding $\overline{\Q} \hookrightarrow \overline{\Qp}$, and thus they have a well-defined valuation. We say that $(\phi_1,...,\phi_h)$ (resp. $\Phi$) has \emph{small slope} if $v(a_{\pri}) < (k+1)/e_{\pri}$ for all $\pri|p$, where $e_{\pri}$ is the ramification index of $\pri$ in $K$. Note that this is precisely the condition that allows us to lift $\phi$ using one of the control theorems above. We say that the \emph{slope} is $(v(a_{\pri}))_{\pri|p}$.
\end{mdef}
Let $\ff$ be an ideal of $\roi_{K}$ with $\ff |(p^\infty)$, and for each $b \newmod{\ff}$, take an element $d_b \in \roi_K$ such that $d_b \in I_1,...,I_h$ and $d_b \equiv b \newmod{\ff}$ using the Chinese Remainder Theorem. Then consider the locally polynomial function
\[P_{b,\ff}^{q,r}(z) \defeq z^q\overline{z}^r\mathbbm{1}_{b\newmod{\ff}},\]
on $\roi_{K,p}$, where $\mathbbm{1}_{b\newmod{\ff}}$ is the indicator function for the minimal open subset of $\roi_{K,p}^\times$ containing the image of $b + \ff \subset \roi_K^\times$ under the canonical embedding of $\roi_K^\times$ into $\roi_{K,p}^\times$. Finally, for each $i$, we write $\ff I_i = (\alpha_i)I_{j_i}.$\\
\\
We define an operator $U_{\ff}$ as follows:
\[U_{\ff} \defeq \prod\limits_{\pri^n||\ff}U_{\pri}^n .\]
Take some small slope classical Bianchi eigensymbol $(\phi_1,...,\phi_j) \in \symb_{\Omega_1(\n)}(V_{k,k}^*(L))$ with $U_{\ff}$-eigenvalue $\lambda_{\ff}$, and lift it to an overconvergent eigensymbol $(\Psi_1,...,\Psi_h)$ using the control theorem. Evaluating $(\Psi_1,...,\Psi_h)$ at $\{0\}-\{\infty\}$, and then evaluating the resulting distribution at this polynomial, we obtain
\begin{align}\label{overconvergentvalue}(\Psi_1,..,\Psi_h)&(\zeroinf)(P_{b,\ff}^{q,r}) = \lambda_{\ff}^{-1}[(\Psi_1,...,\Psi_h)|U_{\ff}](\zeroinf)(z^q\overline{z}^r\mathbbm{1}_{b\newmod{\ff}})\notag\\
&= \lambda_{\ff}^{-1}\left(\sum\limits_{b\newmod{\ff}}\Psi_{j_1}\bigg|\matrd{1}{d_b}{0}{\alpha_1},...,\sum\limits_{b \newmod{\ff}}\Psi_{j_h}\bigg|\matrd{1}{d_b}{0}{\alpha_h}\right)\zeroinf\notag \\
&= \lambda_{\ff}^{-1}\bigg(\Psi_{j_i}\left(\{d_b/\alpha_i\}-\{\infty\}\right)[(\alpha_i z + d_b)^q(\overline{\alpha_i z} + \overline{d_b})^r]\bigg)_{i=1}^h,
\end{align}
where here the sum for $U_{\ff}$ is `absorbed' by the indicator function; indeed, we have
\begin{align*}
\left[\matrd{1}{d_b}{0}{\sigma},\matrd{1}{\overline{\partial}}{0}{\overline{\sigma}}\right]\cdot\left(z^q\overline{z}^r\mathbbm{1}_{b\newmod{\ff}}\right)(x,y) &= z^q\overline{z}^r\mathbbm{1}_{b\newmod{\ff}}(\partial+\sigma x,\overline{\partial}+\sigma y)\\
&=0 \text{ unless }\partial\equiv b\newmod{\ff}.
\end{align*}
Suppose now that for $d \in K$ and $\alpha \in K^\times$ we set
\begin{equation}\label{defncij}
\phi(\{d/\alpha\}-\{\infty\}) = \sum\limits_{i,j=0}^kc_{i,j}\left(\frac{d}{\alpha}\right)\left(\Yc-\frac{d}{\alpha}\Xc\right)^{k-i}\Xc^i\left(\Ycbar-\frac{\overline{d}}{\overline{\alpha}}\Xcbar\right)^{k-j}\Xcbar^j,\end{equation}
where $\Xc^i\Yc^{k-i}\Xcbar^j\Ycbar^{k-j}$ is the basis element of $V_{k,k}^*(L)$ such that 
\[\Xc^i\Yc^{k-i}\Xcbar^j\Ycbar^{k-j}(X^IY^{k-I}\Xbar^J\Ybar^{k-J}) = \delta_{iI}\delta_{jJ}.\]
Note that this is chosen so that under the change of basis for $V_{k,k}(L)$ defined by
\[X^iY^{k-i}\Xbar^j\Ybar^{k-j} \longmapsto (\alpha X + dY)^iY^{k-i}(\overline{\alpha}\Xbar+\overline{d}\Ybar)^j\Ybar^{k-j},\]
the corresponding change of dual basis is given by
\[\Xc^i\Yc^{k-i}\Xcbar^j\Ycbar^{k-j} \longmapsto \Xc^i(\alpha\Yc - d\Xc)^{k-i}\Xcbar^j(\overline{\alpha}\Ycbar - \overline{d}\Xcbar)^{k-j}.\]
When we substitute equation (\ref{defncij}) into (\ref{overconvergentvalue}), noting that this makes sense by definition of $\Psi$ as a lift of $\phi$, and using the obvious dictionary between the two spaces $\mathcal{D}_{k,k}(L)$ and $V_{k,k}^*(L)$, we find that
\begin{align*}\Psi_{i}(\zeroinf)(P_{b,\ff}^{q,r}) &= \lambda_{\ff}^{-1}\alpha_i^{q}\overline{\alpha_i}^{r}c_{q,r}^{j_i}\left(\frac{d_b}{\alpha_i}\right)\\
& = \lambda_{\ff}^{-1}\psi(t_i)^{-1}\psi(t_{j_i})\psi_{\ff}(d_b)^{-1}\psi_{\ff}\left(\frac{d_b}{\alpha_i}\right)c_{q,r}^{j_i}\left(\frac{d_b}{\alpha_i}\right),
\end{align*}
where here we've used that
\[\alpha_i^q\overline{\alpha_i}^r = \psi_\infty(\alpha_i) = \psi(t_i)^{-1}\psi(t_{j_i})\psi_{\ff}(\alpha_i)^{-1}.\]
Note that $\alpha_i^{-1} \in \ff^{-1}I_i^{-1}I_{j_i} \subset \ff^{-1}I_i^{-1}$, so in particular as $b$ ranges over all classes of  $(\roi_K/\ff)^\times$ and as $d_b \in I_i$, we see that $d_b/\alpha_i$ ranges over a full set of coset representatives $[a]$ for $\ff^{-1}/\roi_K$ with $(a)\ff$ coprime to $\ff$. (Note this relies on the fact that we are taking invertible elements $\newmod{\ff}$. Accordingly, the ideal $(d_b/\alpha_i) = \ff^{-1}J$, where $J$ is coprime to $\ff$, and hence $d_b/\alpha_i \notin \roi_K$. It is clear that if $b \neq b' \newmod{\ff}$ then $d_b/\alpha_i$ and $d_{b'}.\alpha_i$ define different classes in $\ff^{-1}/\roi_K)$. Thus the values $c_{q,r}^{j_i}(d_b/\alpha_i)$ are precisely what we need to access the $L$-values, since they occur in an integral formula for the critical values of the $L$-function.

\subsection{Ray Class Groups}
The $p$-adic $L$-function of a modular form should be a function on characters in a suitable sense. To make this more precise, we introduce ray class groups.
\begin{mdef}
Let $K$ be a number field with ring of integers $\roi_K$, and take an ideal $\ff \subset \roi_K$. The \emph{ray class group of $K$ modulo $\ff$}, denoted $\cl(K,\ff)$, is the group $I_{\ff}$ of fractional ideals of $K$ that are coprime to $\ff$ modulo the group $K_{\ff}^1$ of principal ideals that have a generator congruent to 1 mod $\ff$. We can also define the ray class group adelically; if we let $U(\ff) = 1 + \ff \widehat{\roi_F},$ then
\[\cl(K,\ff) \cong K^\times\backslash\A_K^\times/U(\ff)\C^\times.\]
\end{mdef}
The ray class group fits into a useful exact sequence; we have
\[0\longrightarrow \roi_K^\times \newmod{\ff} \longrightarrow (\roi_K/\ff)^\times \labelrightarrow{\beta} \cl(K,\ff) \longrightarrow \cl(K) \longrightarrow 0,\]
where the map $\beta$ takes an element $\alpha + \ff$ to the class $(\alpha) + K_{\ff}^1$ and the surjection is the natural quotient map. Now let $\ff|(p^\infty)$. Piecing this together as we let $n$ vary, taking the inverse limit of this family of exact sequences, we obtain an exact sequence
\[0 \longrightarrow \roi_K^\times \longrightarrow (\roi_K\otimes_{\Z}\Z_p)^\times \longrightarrow \cl(K,p^\infty) \labelrightarrow{\delta} \cl(K) \longrightarrow 0,\]
where here 
\[\cl(K,p^\infty) \defeq \lim\limits_{\longleftarrow}\cl(K,p^\infty)\]
is defined to be the inverse limit.\\

\subsection{Grossencharacters as Characters of $\cl(K,p^\infty)$}
Let $\mathfrak{g}$ be some ideal of $\roi_K$ that is coprime to $(p)$. The following theorem is due to Weil; note here that we consider $\cl(K,\mathfrak{g}p^\infty)$ as a quotient of the adeles in the natural way.
\begin{mthm}[Weil] \label{weilassoc} There is a bijection between algebraic Grossencharacters of conductor dividing $\mathfrak{g}p^\infty$ and locally algebraic characters of $\cl(K,\mathfrak{g}p^\infty)$ such that if $\psi$ corresponds to $\psi_{p-\mathrm{fin}}$, then
\[\psi\bigg|_{\left(\A_K^{(p,\infty)}\right)^\times} = \psi_{p-\mathrm{fin}}\bigg|_{\left(\A_K^{(p,\infty)}\right)^\times},\]
where here consider the adelic definition of the ray class group and we restrict to the adeles away from the infinite place and the primes above $p$.
\end{mthm}
The correspondence is simple to describe. If $\psi$ is an algebraic Grossencharacter of conductor $\ff|\mathfrak{g}p^\infty$ and infinity type $(q,r)$, then to $\psi$ we associate a $K^\times$-invariant function
\[\psi_{p-\text{fin}}:(\A_K^\times)_f \longrightarrow \C_p^\times\]
by fixing an isomorphism $\C\cong\C_p$ and defining
\[\psi_{p-\text{fin}}(x) \defeq \psi_f(x)\sigma_p(x),\]
where
\[\sigma_p(x) \defeq \left\{\begin{array}{ll}x_{\pri}^q x_{\pribar}^r &: p \text{ splits as }\pri\pribar,\\
x_{\pri}^q \overline{x_{\pri}^r} &: p\text{ inert or ramified}.\end{array}\right.\]
It is simple to check that this is trivial on $K^\times$. By unravelling the definitions, we see that for $\alpha \in K^\times$, and an idele $x_{\alpha,p}$ defined by
\[(x_{\alpha,p})_{\pri} \defeq \left\{\begin{array}{ll}\alpha&: \pri|(p)\\ 1&:\text{otherwise,}\end{array}\right.\]
we have
\[\psi_{p-\text{fin}}(x_{\alpha,p}) = (\psi_{p-\text{fin}})_{(p)}(\alpha) = \psi_{(p)}(\alpha)\alpha^q\overline{\alpha}^r.\]
\\
This association now gives the construction of Theorem \ref{weilassoc}.

\subsection{Constructing the $p$-adic $L$-function}
The $p$-adic $L$-function of a modular form over a field $K$ is naturally a locally analytic distribution on $\cl(K,p^\infty)$. We've reformulated the complex $L$-function as a function of Grossencharacters; we'll now prove that we can construct such a locally analytic distribution $\mu_p$ that satisfies the interpolation property that, for any Grossencharacter of conductor $\ff|(p^\infty)$ and infinity type $0 \leq (q,r) \leq (k,k)$, we have
\[\mu_p(\psi_{p-\mathrm{fin}}) = \left[\frac{(-1)^{k+q+r}2\psi_{\ff}(x_{\ff})\lambda_{\ff}}{\psi(x_{\ff})Dw\tau(\psi^{-1})}\right]^{-1}\Lambda(\Phi,\psi),\]
where $\lambda_{\ff}$ is the $U_{\ff}$-eigenvalue of $\Phi$. If $\Phi$ has small slope, then in conjunction with the admissibility condition, this defines $\mu_p$ uniquely (see, for example, \cite{Loe14} for this result in the weight (0,0) case). The \emph{$p$-adic $L$-function of $\Phi$} will then be defined to be the distribution $\mu_p$. To make this tie in with the work of Section \ref{evalatzeroinf}, we describe the character $\psi_p$ explicitly as a locally analytic function on the ray class group.\\
\\
First, we describe what $\cl(K,p^\infty)$ actually looks like. By choosing a set of representatives for the class group, we are choosing a section of the map $\cl(K,p^\infty) \rightarrow \cl(K)$, and thus, going back to the exact sequence above, we can identify $\cl(K,p^\infty)$ with a disjoint union of $h$ copies of $\roi_{K,p}^\times/\roi_K^\times$, indexed by our class group representatives. On each of the $h$ components, the character $\psi_{p-\mathrm{fin}}$ gives a locally polynomial function on $\roi_{K,p}^\times$; this is given by
\begin{align*}P_i(z) &= \psi(t_i)\sum\limits_{b \in (\roi_K/\ff)^\times}\psi_{\ff}(b)\eta(z)\mathbbm{1}_{b\newmod{\ff}}\\
& = \psi(t_i)\sum\limits_{b \in (\roi_K/\ff)^\times}\psi_{\ff}(b)P_{a,\ff}^{q,r}(z),
\end{align*}
where $\eta$ is the character sending $z \mapsto z^q\overline{z}^r$, and $\mathbbm{1}_{b\newmod{\ff}}$ is the indicator function for the minimal open subset of $\roi_{K,p}^\times$ containing the image of $b \newmod{\ff} \subset \roi_K^\times$ under the canonical embedding of $\roi_K^\times$ into $\roi_{K,p}^\times$.\\
\\
With the above data, combined with equation (\ref{finallformula}), we are in a position to construct our distribution on $\cl(K,p^\infty).$ We recall the construction. Let $\n$ be an ideal divisible by $(p)$, and let $\Phi$ be a cuspidal Bianchi eigenform of level $\Omega_1(\n)$ and weight $(k,k)$. Then to such a function we associate a modular symbol $(\phi_1, ..., \phi_h)$, where $h$ is the class number. By scaling by the periods of the eigenform, we can consider each $\phi_i$ as an element of Hom$_{\Gamma_i}(\Delta_0,V_{k,k}(L))$ for some extension $L/\Q_p$. In the case where $\Phi$ (and hence each $\phi_i$) has small slope $(h_{\pri})_{\pri|p}$, we can lift these symbols to give an overconvergent modular symbol $(\Psi_1,...,\Psi_h)$ in $\oplus_{i=1}^h$Hom$_{\Gamma_i}(\Delta_0,\mathcal{D}_{k,k}(L))$. The values taken by such symbols are then $(h_{\pri})_{\pri|p}$-admissible in the sense of the previous section.\\
\\
Define, for each $i$, a distribution
\[\mu_i \defeq \Psi_i(\{0\}-\{\infty\})\bigg|_{\roi_{K,p}^\times}.\]
Then we know that, for $\ff|(p^\infty),$ 
\[\mu_i\left(P_{b,\ff}^{q,r}(z)\right) = \lambda_{\ff}\psi(t_i)\psi(t_{j_i})\psi_{\ff}(d_b)^{-1}\psi_{\ff}(d_b/\alpha_i)c^{j_i}_{q,r}(d_b/\alpha_i),\]
where $\lambda_{\ff}$ is the $U_{\ff}$-eigenvalue of $\Phi$, and which leads us to define a distribution 
\[\mu_p \defeq \sum\limits_{i=1}^h\mu_i \mathbbm{1}_i\]
on $\cl(K,p^\infty)$, where here $\mathbbm{1}_i$ means the indicator function for the component of $\cl(K,p^\infty)$ corresponding to $I_i$. Then take a Grossencharacter $\psi$ of infinity type $(q,r)$ -- where $0 \leq q,r \leq k$ -- and conductor $\ff|(p^\infty)$. Then we have
\begin{align*}\mu_p(\psi_{p-\mathrm{fin}}) &= \mu_p\left(\sum\limits_{i=1}^hP_i(z)\mathbbm{1}_{i}\right) = \sum\limits_{i=1}^h\sum\limits_{b \in (\roi_K/\ff)^\times}\psi(t_i)\psi_{\ff}(d_b)\mu_i\left(P_{b,\ff}^{q,r}(z)\right)\\
& = \lambda_{\ff}^{-1}\sum\limits_{i=1}^h\psi(t_{j_i})\sum\limits_{b \in (\roi_K/\ff)^\times}\psi_{\ff}\left(\frac{d_b}{\alpha_i}\right)c_{q,r}^{j_i}\left(\frac{d_b}{\alpha_i}\right)\\
& = \left[\frac{(-1)^{k+q+r}2\psi_{\ff}(x_{\ff})\lambda_{\ff}}{\psi(x_{\ff})Dw\tau(\psi^{-1})}\right]^{-1}\Lambda(\Phi,\psi),
\end{align*}
using equation (\ref{finallformula}). In this equation, recall that $\psi_{\ff} = \prod_{\pri|\ff}\psi_{\pri},$ the idele $x_{\ff}$ is as defined in Section \ref{lfnchar}, $\lambda_{\ff}$ is the $U_{\ff}$-eigenvalue of $\Phi$, $-D$ is the discriminant of $K$, $w$ is the size of the unit group of $K$ and $\tau(\psi^{-1})$ is a Gauss sum as defined in Section \ref{lfunction}. This interpolation property, in addition to the admissibility condition described above, means we can now prove:
\begin{mthm}\label{padiclfn}
Let $K/\Q$ be an imaginary quadratic field of class number $h$ and discriminant $-D$, and let $p$ be a rational prime. Let $\Phi$ be a cuspidal Bianchi eigenform of weight $(k,k)$ and level $\Omega_1(\n)$, where $(p)|\n$, with $U_\pri$-eigenvalues $a_\pri$, where $v(a_{\pri})<(k+1)/e_{\pri}$ for all $\pri|p$. Then there exists a locally analytic distribution $\mu_p$ on $\cl(K, p^\infty)$ such that for any Grossencharacter of $K$ of conductor $\ff|(p^\infty)$ and infinity type $0 \leq (q,r) \leq (k,k)$, we have
\[\mu_p(\psi_{p-\mathrm{fin}}) = \left[\frac{(-1)^{k+q+r}2\psi_{\ff}(x_{\ff})\lambda_{\ff}}{\psi(x_{\ff})Dw\tau(\psi^{-1})}\right]^{-1}\Lambda(\Phi,\psi),\]
where $\psi_{p-\mathrm{fin}}$ is as in Theorem \ref{weilassoc}. The distribution $\mu_p$ is $(h_{\pri})_{\pri|p}$-admissible in the sense of Definitions \ref{admissible} and \ref{splitadmissible}, where $h_{\pri} = v_p(a_{\pri})$, and hence is unique.\\
\\
We call $\mu_p$ the \emph{$p$-adic $L$-function of $\Phi$}.
\end{mthm}
\begin{proof}
The eigenform $\Phi$ corresponds to a collection of $h$ cusp forms $\f^1,...,\f^h$ on $\uhs$; associate to each $\f^i$ a classical Bianchi eigensymbol $\phi_{\f^i}$ with coefficients in a $p$-adic field $L$, and lift each to its corresponding unique overconvergent Bianchi eigensymbol $\Psi_{i}$. Define 
\[\mu_i \defeq \Psi_i(\{0\}-\{\infty\})\bigg|_{\roi_{K,p}^\times},\] 
and define a locally analytic distribution $\mu_p$ on $\cl(K,p^\infty)$ by $\mu_p \defeq \sum_{i=1}^h\mu_i\mathbbm{1}_i.$ Then by the work above, $\mu_p$ satisfies the interpolation and admissibility properties.
\end{proof}
\begin{mrem}As an example of where this theorem applies, suppose $p$ splits in $K$ and let $E/K$ be a modular elliptic curve with supersingular reduction at both primes above $p$. Then to $E$ we can associate a modular symbol which will have slope $1/2$ at each of the primes above $p$. Accordingly, our construction will give the $p$-adic $L$-function of $E$.
\end{mrem}

\bibliography{references}{}
\bibliographystyle{alpha}
\end{document}